\input amstex
\documentstyle{amsppt}
\magnification=\magstep1

\NoBlackBoxes
\TagsAsMath

\pageheight{9.0truein}
\pagewidth{6.5truein}

\long\def\ignore#1\endignore{\par DIAGRAM\par}
\long\def\ignore#1\endignore{#1}

\ignore
\input xy \xyoption{matrix} \xyoption{arrow}
          \xyoption{curve}  \xyoption{frame}
\def\edge{\ar@{-}}
\def\dttdar{\ar@{.>}}
\def\drbl{\save+<0ex,-2ex> \drop{\bullet} \restore}

\def\dashedge{\ar@{--}}

\def\dshdar{\ar@{-->}}
\endignore

\def\la{{\Lambda}}
\def\lamod{\Lambda\text{-}\roman{mod}}

\def\Lamod{\Lambda\text{-}\roman{Mod}}
\def\Amod{A\text{-}\roman{Mod}}
\def \len{\operatorname{length}}

\def\SS{{\Bbb S}}

\def\hom{\operatorname{Hom}}

\def\soc{\operatorname{soc}}

\def\pdim{\operatorname{proj\,dim}}

\def\lfindim{\operatorname{l.findim}}
\def\Lfindim{\operatorname{l.Findim}}
\def\rfindim{\operatorname{r.findim}}
\def\Rfindim{\operatorname{r.Findim}}

\def\Ker{\operatorname{Ker}}

\def\Add{\operatorname{Add}}
\def\Ext{\operatorname{Ext}}

\def\End{\operatorname{End}}

\def\can{\operatorname{can}}

\def\eps{\epsilon}

\def\pdim{\operatorname{p\, dim}}
\def\gldim{\operatorname{gl\, dim}}

\def\A{{\Cal A}}
\def\B{{\Cal B}}
\def\C{{\Cal C}}

\def\P{{\Cal P}}

\def\S{{\sigma}}

\def\eps{{\varepsilon}}

\def\etilde{\widetilde{e}}
\def\Itilde{\widetilde{I}}
\def\Jtilde{\widetilde{J}}
\def\Stilde{\widetilde{S}}
\def\Qtilde{\widetilde{Q}}
\def\latilde{\widetilde{\Lambda}}
\def\Etilde{\widetilde{E}}
\def\Mtilde{\widetilde{M}}
\def\epstilde{\widetilde{\epsilon}}

\def\Mod{\operatorname{Mod}}
\def\mod{\operatorname{mod}}

\def\pinf{\operatorname{\P^{< \infty}}}
\def\pinflamod{\operatorname{\P^{< \infty}(\lamod)}}
\def\Pinflamod{\operatorname{\P^{< \infty}(\Lamod)}}
\def\pinfamod{\operatorname{\P^{< \infty}(A\text{-}mod)}}
\def\Pinfamod{\operatorname{\P^{< \infty}(A\text{-}Mod)}}

\def\pinfmodb{\operatorname{\P^{< \infty}(mod\text{-}B)}}
\def\pinflatilde{\operatorname{\P^{<
\infty}(mod\text{-}\latilde)}}
\def\Pinflatilde{\operatorname{\P^{<
\infty}(Mod\text{-}\latilde)}}
\def\add{\operatorname{add}}
\def\filt{\operatorname{filt}}
\def\Filt{\operatorname{Filt}}

\def\Hom{\operatorname{Hom}}

\def\id{\operatorname{id}}

\def\add{\operatorname{add}}
\def\Add{\operatorname{Add}}

\def\fcog{\operatorname{fcog}}
\def\Coker{\operatorname{Coker}}

\def\AgHaLuUn{{\bf 1}}
\def\AnTr{{\bf 2}}
\def\AuGr{{\bf 3}}
\def\AuRe{{\bf 4}}
\def\AuSm{{\bf 5}}
\def\BHT{{\bf 6}}
\def\BHZ{{\bf 7}}
\def\CPS{{\bf 8}}
\def\Dla{{\bf 9}}
\def\DH{{\bf 10}}
\def\DHL{{\bf 11}}
\def\Fri{{\bf 12}}
\def\Hap{{\bf 13}}
\def\HaUn{{\bf 14}}
\def\dom{{\bf 15}}
\def\menace{{\bf 16}}
\def\HuSm{{\bf 17}}
\def\Miya{{\bf 18}}
\def\PaXi{{\bf 19}}
\def\Rin{{\bf 20}}
\def\Sma{{\bf 21}}

\topmatter

\title Strongly tilting truncated path algebras
\endtitle

\author A. Dugas and B. Huisgen-Zimmermann
\endauthor

\address Department of Mathematics, University of California, Santa
Barbara, CA 93106-3080 \endaddress

\thanks The  second author was partly supported by a grant from the
National Science Foundation.
\endthanks

\abstract For any truncated path algebra $\la$, we give a structural
description of the modules in the categories $\pinflamod$ and
$\Pinflamod$, consisting of the finitely generated (resp. arbitrary)
$\la$-modules of finite projective dimension.  We deduce that these
categories are contravariantly finite in $\lamod$ and $\Lamod$,
respectively, and determine the corresponding minimal
$\pinf$-approximation of an arbitrary $\la$-module from a projective
presentation.  In particular, we explicitly construct   --  based on
the  Gabriel quiver
$Q$ and the Loewy length of $\la$  -- the basic strong tilting module
$_\la T$  (in the sense of Auslander and Reiten) which is coupled with
$\pinflamod$ in the contravariantly finite case.  A main topic is the
study of the homological properties of the corresponding tilted algebra
$\latilde = \End_\la(T)^{\text{op}}$, such as its finitistic
dimensions and the structure of its modules of finite projective
dimension.  In particular, we characterize, in terms of a
straightforward condition on
$Q$, the situation where the tilting module
$T_{\latilde}$ is strong over $\latilde$ as well.  In this
$\la$-$\latilde$-symmetric situation, we obtain sharp results on the
submodule lattices of the objects in
$\pinf(\Mod$-$\latilde)$, among them  a certain heredity property; it
entails that any module in
$\pinf(\Mod$-$\latilde)$ is an extension of a projective module by a
module all of whose simple composition factors belong to
$\pinf(\mod$-$\latilde)$.
\endabstract

\endtopmatter

\document

\head 1. Introduction and terminology \endhead

We let $\la = KQ/I$ be a truncated path algebra of Loewy length
$L+1$ for some positive integer $L$, meaning that $KQ$ is the path
algebra of a quiver
$Q$ with coefficients in a field $K$ and $I \subseteq KQ$ the ideal
generated by all paths of length
$L+1$.  Provided that $K$ is algebraically closed, the class of
truncated path algebras includes all basic hereditary algebras, as well
as all basic algebras with vanishing radical square.  Since we place no
restrictions beyond finiteness on the quiver $Q$, algebraic closedness
of the base field moreover entails that every finite dimensional
$K$-algebra is Morita equivalent to a factor algebra of a truncated path
algebra.  Our results do not require any hypothesis on $K$, however. By
$\lamod$ (resp\.
$\Lamod$), we denote the category of all finitely generated (resp\. all)
left
$\la$-modules.

In Section 3, we structurally characterize the objects in the
subcategories
$\pinflamod$ and $\Pinflamod$, consisting of the modules of finite
projective dimension in the categories $\lamod$ and $\Lamod$,
respectively.  Our description rests on the following two facts: Every
$\la$-module
$M$ contains a unique largest submodule $U(M)$ all of whose composition
factors have finite projective dimension.  Moreover, there
are finitely many local $\la$-modules $\A_i$ giving rise to the
following test for finiteness of the projective dimension:  Namely,
$M$ belongs to $\Pinflamod$ if and only if $M/U(M)$ is a direct sum of
copies of the $\A_i$; see Theorem 3.1 and Corollary 3.3 for more
precision.  As one byproduct of this result, we see that the category
$\Pinflamod$ is closed under top-stable submodules (a module $N
\subseteq M$ is a {\it top-stable submodule of $M$\/} in case $JN = JM
\cap N$, where
$J$ is the Jacobson radical of $\la$).

As another consequence of the mentioned ``homological subdivision" of
$\la$-modules, we find that the categories $\pinflamod$ and
$\Pinflamod$  are contravariantly finite in $\lamod$ and $\Lamod$,
respectively (Theorems 4.1 and 4.2).  The $\A_i$ mentioned above turn
out to be the minimal
$\pinflamod$-approximations of the simple modules of infinite projective
dimension; their structure is immediate from the quiver and Loewy length of
$\la$.   This adds another instance to the short list of known
classes of finite dimensional algebras
$A$ whose categories $\pinfamod$ are consistently contravariantly
finite:  So far, this has been established whenever $A$ is stably
equivalent to a hereditary algebra (see
\cite{\AuRe, p\. 130}), or else when $A$ is left serial (see
\cite{\BHZ}); the former class, in turn, contains the radical-square
zero algebras.   Over a truncated path algebra
$\la$, the minimal
$\Pinflamod$-approximation of any $\la$-module $M$ is readily accessible
from a minimal projective presentation of $M$; the connection is
described in Theorem 4.2.  We conclude Section 4 by showing that
contravariant finiteness of
$\pinflamod$ can alternatively be derived from a theorem of Smal\o\ in
\cite{\Sma}.

The homological picture of $\la$, which started emerging in
\cite{\DHL}, where the homological dimensions of
$\la$ were pinned down in terms of the quiver of the algebra and
of its Loewy length, is based on a bicoloring
of the vertices of the quiver $Q$, precyclic versus non-precyclic.  We call a vertex
$e$ {\it precyclic\/} if there is an oriented path which starts in $e$
and ends on an oriented cycle.  It is easily seen that a vertex $e_i$ of
$Q$ is precyclic if and only if the corresponding simple left module
$S_i =
\la e_i / Je_i$ has infinite projective dimension (for a more general
result on projective dimensions of local $\la$-modules in terms of their
tops, see \cite{\DHL, Theorem 2.6}).  The bicoloring continues to be
pivotal in our description of the unique basic tilting module $T$ of
$\la$ which is  Ext-injective in
$\pinflamod$.  In \cite{\AuRe, Section 6}, such a tilting module was
proved to exist, over a finite dimensional algebra $A$ say, precisely
when $\pinfamod$ is contravariantly finite; in case of existence, it was
dubbed the {\it strong\/} ({\it basic\/}) {\it tilting module\/} in
$A$-$\mod$ -- see the beginning of Section 5 for background.  In our
scenario, that is, over a truncated path algebra $\Lambda$, the
structure of the strong tilting module $_\la T$ can be pinned down (see
Theorem 5.3).  In particular, our characterization permits us to
construct $T$ from the basic data, $Q$ and
$L$; for concrete illustrations we refer to Examples 5.6.
 As a consequence, the quiver and relations of the tilted algebra
$\latilde := \End_\la(T)^{\text{op}}$ can in turn be determined from
these data, albeit with some computational effort; instead of giving a
cumbersome formal algorithm, we include two examples at the end (Section
9).

Typically, the tilted algebra $\latilde$ has higher Loewy length than
$\la$, and the basic oriented cycles of its quiver $\Qtilde$ may increase
in number and length when compared with those of
$Q$.  In Section 9, we present an example of a truncated path algebra
$\la$ with Loewy length $3$ and a quiver having only one basic oriented
cycle, while
$\latilde$ has Loewy length $7$ and four distinct basic oriented
cycles.  Moreover:  Whereas path length in $Q$ clearly induces a grading
of $\la$, the analogue for $\Qtilde$ and $\latilde$ fails, in general.
However, an in-depth study of
$\latilde$ does reveal a natural grading that stems from a valuation of
$\Qtilde$ in general, a point which will be solidified in a sequel to
this article.

In Section 6, we describe filtrations of the objects in the categories
$_\la T^\perp \subseteq \lamod$ and
$^{\perp} (_{\latilde}DT) \subseteq \latilde$-$\mod$, finding parallels
with the theory of quasi-hereditary algebras.  As is the case for the
latter algebras, any truncated path algebra $\la$ is standardly
stratified (in the weak sense of Cline, Parshall and Scott
\cite{\CPS}), relative to a suitable pre-order on the set of simples.
While such stratifications are much coarser than those introduced by
Dlab
\cite{\Dla} under the same name, Frisk has recently shown that many
results known for Dlab's standardly stratified algebras can be extended
to the more general situation
\cite{\Fri}.  We illustrate this theory in the case of a
truncated path algebra (see Theorem 6.1 and Remark 6.2) and refer
to Remarks 3.4 and 8.6 for more information on the connection.

Next, we proceed to a structural exploration of the objects in
$\pinf(\Mod$-$\latilde)$.  The most symmetric and transparent situation
occurs when
$Q$ has no precyclic source.  In Section 7, we show this  condition to
be equivalent to the requirement that the tilting bimodule $_\la
T_{\latilde}$ be strong on both sides (Theorem 7.2).  Thus, the algebras
$\latilde$ obtained by strongly tilting truncated path algebras with
quivers devoid of precyclic sources constitute yet another class of
algebras $A$ enjoying contravariant finiteness of $\pinf(\mod
\text{-}A)$.  On one hand, these algebras $\latilde$ are considerably
more complex in structure than the aforementioned examples.  On the
other hand, further applications of tilting theory yield substantial
information on the objects of
$\pinf(\Mod$-$\latilde)$.  The sharpest structural results can be
found in Theorems 8.2 and 8.5. We remark that the results of Sections
7--9 do not depend on Section 6, but are
linked directly to Sections 3--5.
\bigskip

\noindent {\it Terminology.}  Let $K$ be an arbitrary field,
$Q$ a finite quiver and $L$ an integer $\ge 1$. {\it Throughout, $\la =
KQ/I$ will stand for a truncated path algebra\/} with radical $J$ and
$J^{L+1} = 0$; in other words, $I
\subseteq KQ$ will denote the ideal generated by all paths of length
$L+1$ in
$Q$.  (Whenever we address algebras that are not necessarily of this
type, we will use a different notation.)  The set of vertices of
$Q$ will be identified with the paths of length zero in
$KQ$, and further with a full sequence $e_1, \dots, e_n$ of primitive
idempotents of
$\la$.  Our convention for multiplying paths $p, q \in KQ$ is as follows:
$pq$ stands for ``$p$ after $q$".  In keeping with this convention, we
call a path $p'$ an {\it initial\/} (or {\it terminal\/}) {\it
subpath\/} of a path $p$ if $p = p'' p'$ (or
$p = p' p''$).  Moreover, a {\it path in $\la$\/} is any residue class
$p + I$, where $p$ is a path in $KQ \setminus I$.  It is clearly
unambiguous to carry over the notions of length, starting point and end
point from paths in $KQ$ to paths in $\la$, since the ideal $I$ is
homogeneous with respect to the path-length grading of $KQ$.
Representatives of the simples in
$\lamod$ are
$S_i = \la e_i / J e_i$, $1 \le i \le n$.

A vertex $e_i$ of $Q$ is called {\it precyclic\/} in case there exists a
path in $Q$ which starts in $e_i$ and terminates in a vertex lying on an
oriented cycle.  Dually, $e_i$ is {\it postcyclic\/} if $e_i$ is the
endpoint of a path in $Q$ that starts on an oriented cycle.
Correspondingly, we also refer to the simple module $S_i$ as
{\it precyclic\/} or {\it postcyclic\/}.

An auxiliary concept we use is that of a {\it sequence of top
elements\/} of
$M \in \Lamod$:  We call an element $m \in M \setminus JM$ a {\it top
element\/} of $M$ if it is normed by one of the primitive idempotents,
i.e., $e_i m = m$ for some $i \in \{1,
\dots, n\}$.  A family
$(m_r)_{r\in R}$ of top elements will be called a {\it sequence of top
elements of $M$\/} provided that the residue classes
$m_r + JM$ form a basis for $M/JM$ over $K$.  For instance
$e_1, \dots, e_n$ is a sequence of top elements of the left regular
module $\la$.

By a {\it tilting module\/} $_A T$ over any finite dimensional algebra
$A$ we mean a finitely generated module of finite projective dimension
with
$\Ext^i_A(T, T) = 0$ for all $i
\ge 1$ such that the regular left $A$-module $A$ has a finite
(exact) coresolution
$$0 \rightarrow {}_A A \rightarrow M_0 \rightarrow \cdots
\rightarrow M_t \rightarrow 0$$ 
with $M_i \in \add(T)$ for all $i \ge 0$.

For any subcategory $\Cal{C}$ of $A$-$\mod$, we define its left and right perpendicular subcategories
to be
$${}^{\perp}\Cal{C} = \{M \in A
\text{-}\mod\ |\ \Ext^i_A(M,C) = 0\ \text{for all}\ C \in
\Cal{C} \text{\ and\ } i \ge 1\}$$  and
$$\Cal{C}^{\perp} = \{ M
\in A \text{-}\mod\ |\ \Ext^i_A(C,M) = 0\ \text{for all}\ C \in
\Cal{C} \text{\ and\ } i \ge 1\},$$    respectively.  As usual, given a
module $C$, we write $C^{\perp}$ for $\{C\}^{\perp}$.

Finally, we call the invariants
$$\lfindim A = \sup\{\pdim M \mid M \in \pinfamod\}$$ and
$$\Lfindim A = \sup\{\pdim M \mid M \in \Pinfamod\}$$ the left {\it
little\/} and {\it big finitistic dimensions\/} of $A$, respectively.

\head 2. Prerequisites \endhead

Let $A$ be any finite dimensional algebra and $\C$ a subcategory of the
category $A$-$\mod$ of finitely generated left $A$-modules which is
closed under direct summands.  Recall that, according to Auslander and
Smal\o\
\cite{\AuSm}, the subcategory $\C$ is said to be contravariantly finite
in $A$-$\mod$ in case every object
$M$ in $A$-$\mod$ has a (right) $\C$-approximation:  Such an
approximation is a  morphism $\psi: B \rightarrow M$ with $B$ in
$\C$, such that every homomorphism $C \rightarrow M$ with $C$ in $\C$
factors through $\psi$.  By a slight abuse of language, one also refers
to the object $B$ as a $\C$-approximation of $M$ in this case.
Whenever
$M$ has a (right) $\pinfamod$-approximation, there is a ``best", that
is, {\it minimal\/} one, say $\phi: B(M) \rightarrow M$; it is
characterized by the property that, for any endomorphism $g$ of $B(M)$,
the equality
$\phi \circ g = \phi$ forces $g$ to be an automorphism of
$B(M)$.  It is easily checked that $B(M)$ is isomorphic to a direct
summand of any $\C$-approximation of $M$; in particular, $B(M)$ has
minimal $K$-dimension among the approximations of $M$.

Subsequently, this terminology was extended to arbitrary summand-closed
subcategories $\C$ of the big module category
$A$-$\Mod$.  For $\C$ to be {\it contravariantly finite\/} in
$A$-$\Mod$, we require that every left $A$-module should have a (right)
$\C$-approximation; such an approximation is defined as above, on
waiving all conditions involving $K$-dimensions.  The definition of
minimal approximations, in turn, follows the above pattern.

We will briefly refer to a ``$\C$-approximation"  when we mean a ``right
$\C$-aproximation".  Since we do not consider any left approximations in
this paper, this will not lead to ambiguities.

Here, we are primarily interested in the special cases $\C =
\pinfamod$ and
$\C = \Pinfamod$.  By \cite{\AuRe}, the
category
$\pinfamod$ is contravariantly finite in $A$-$\mod$ if and only if every
simple left $A$-module has a
$\pinfamod$-approximation.  In the positive case, suppose that
$\A_1, \dots, \A_n$ are the minimal approximations of the simple left
$A$-modules.  The objects in
$\pinfamod$ are then precisely the direct summands of modules that have
(finite) filtrations with consecutive factors in
$\{\A_1, \dots,
\A_n\}$.  As a consequence, $\lfindim A$ equals the maximum of the
projective dimensions of the
$\A_i$.   Due to \cite{\HuSm}, the objects in $\Pinfamod$ are direct
limits of finitely generated modules of finite projective dimension in
case
$\pinfamod$ is contravariantly finite, which yields $\Lfindim A =
\lfindim A$.
\bigskip

{\it We now specialize to the situation where $A = \la$ is a truncated
path algebra as above\/}.  In this context, we recall a useful
homological tool, the skeleton of a
$\la$-module $M$.  (It is defined in general  --  see
\cite{\BHT}  --  but the general definition can be simplified in the
truncated case.)  A skeleton is a special basis, reflecting the
$KQ$-structure of $M$. As we will see, it provides a convenient
coordinate system for the exploration of the structure of $M$ -- of
homological properties in the present context.  In the case of a
truncated path algebra, any skeleton of $M$ completely determines the
syzygies of
$M$ up to isomorphism, for instance; see Theorem 2.3 below.

\definition{ Definition 2.1. Skeleton of a $\Lambda$-module
$M$}  Fix  a projective cover $P$ of $M$, say $P = \bigoplus_{r
\in R} \Lambda z_r$, where each $z_r$ is one of the idemptents in
$\{e_1, \dots, e_n\}$ tagged with a place number $r$.  We will refer to
the family
$(z_r)_{r \in R}$ as {\it the distinguished sequence of top elements of
$P$\/}.  A {\it path of length\/}
$l$ in $P$ is any nonzero element $p z_r \in P$, where $p$ is a path of
length $l$ in $\la$ (see the first paragraph under terminology above;
moreover, note that $pz_r \ne 0$ forces $p$ to start in the vertex
$e(r)$ norming the top element $z_r$ of
$P$, that is, satisfying $e(r) z_r = z_r$).  Given any set
$\sigma$ of paths in $P$, we denote by
$\sigma_l$ the subset consisting of the paths of length $l$ in
$\sigma$.

\roster
\item"(a)"  A set $\sigma$ of paths of length at most $L$ in
$P$ is a {\it skeleton\/}  of $M$ (in $P$), in case there exists an
epimorphism $f:  P \rightarrow M$ such that, for each
$l
\leq L$, the family of residue classes $f(p z_r)  + J^{l+1} M$, where
$p z_r$ traces the paths in $\sigma_l$, is a
$K$-basis for $J^l M/ J^{l+1} M$.  Moreover, we require that
$\sigma$ be closed under initial subpaths, that is, if $p_2 p_1 z_r \in
\S$, then $p_1 z_r$ in $\sigma$.

\item"(b)" A path $q z_r$ in $P \setminus \sigma$ is called
$\sigma$-{\it critical\/} if it is of the form
$\alpha p z_r$, where $\alpha$ is an arrow and $pz_r$ a path in
$\sigma$ (possibly of length zero).
\endroster
\enddefinition

In particular, the definition entails that any skeleton
$\sigma$ of $M$ in $P$ contains the distinguished sequence of top elements
of $P$ (as the subset $\sigma_0$).  We will typically identify $M$ with a
quotient $P/C$, where $C \subseteq JP$, and focus on subsets $\S \subseteq
P$ which are skeleta of $M$ with respect to the canonical epimorphism
$P \rightarrow P/C$.  For any such skeleton $\S$, the set of residue
classes
$\{p z_r + C \mid p \in \S\}$ is clearly a basis for $M$, the subsets
$\{p z_r + C \mid p \in \S_l \}$ inducing bases for the radical layers
$J^l M/ J^{l+1}M$ for $l \ge 0$.  Note that
the isomorphism class of $P/C$ as a $\la$-module is completely determined
by the expansion coefficients of the elements $q z_r + C$ relative to the
basis $\{p z_r + C \mid p \in \S\}$, where $q z_r$ runs through
the $\sigma$-critical paths in $P$.

It is easily checked
that every $\Lambda$-module
$M$ has at least one skeleton (in any given projective cover $P$ with
distinguished sequence of top elements).  On the other hand, when $P$ is
finite dimensional, the collection of all skeleta of modules
$P/C$ with $C \subseteq JP$ is clearly finite (provided that the
distinguished sequence of top elements of
$P$ is fixed).
Note that $M=\la$, endowed with the distinguished top elements
$e_1,\dots,e_n$, has precisely one skeleton, namely the set of all
paths in $\la$.

Concerning existence, the following
strengthened observation will be useful in Section 3. The final statement of the upcoming lemma is only relevant when the module $M_1$ fails to be finitely generated.

\proclaim{Lemma 2.2} Let $M_1$, $M_2$ be $\la$-modules, not necessarily
finitely generated, and let
$\sigma''$ be a skeleton of $M_2$ {\rm {(}}in some projective cover
$P_2$ of $M_2${\rm {)}}.  If
$M_2$ is an epimorphic image of $M_1$, then $\sigma''$ can be
supplemented to a skeleton $\sigma = \sigma' \sqcup \sigma''$ of
$M_1$ {\rm {(}}in a suitable projective cover of the form $P_1
\oplus P_2$ of
$M_1${\rm {)}}.

Moreover, given any epimorphism $\pi: M_1
\rightarrow M_2$, the skeleton $\sigma= \sigma' \sqcup \sigma''$ of $M_1$ may be chosen {\rm(}dependent on $\pi${\rm)} in such a fashion that $\sigma'$ is empty precisely when $\pi$ is an isomorphism.
\endproclaim

\demo{Proof}  Without loss of generality, $M_2 = M_1/U$, and
$\pi: M_1 \rightarrow M_2$ is the canonical epimorphism.  Moreover, it
is harmless to start with a projective cover $f: P
\rightarrow M_1$ such that $P = P_1 \oplus P_2$, for a suitable
projective module $P_1$, with the property that
$f_2: = \pi \circ f|_{P_2}: P_2 \rightarrow M_2$ is a projective cover
of $M_2$ satisfying condition $(a)$ of the definition, relative to the
skeleton $\sigma''$ of $M_2$.  In other words, we assume that the
elements $f_2(q) + J^{l+1} M_2$, with $q \in
\sigma''_2$ form a basis for $J^l M_2 / J^{l+1} M_2$.  Say
$(z_r)_{r \in R_1}$ and $(z_r)_{r \in R_2}$ are the distinguished
sequences of top elements of $P_1$ and $P_2$, where
$R_1$ and $R_2$ are disjoint index sets; then the union of all the
$z_r$ is the distinguished sequence of $P$.  Note that
$\sigma''_0 = \{z_r \mid r \in R_2\}$, set $\sigma'_0 = \{z_r
\mid r \in R_1\}$, and define $\sigma_0 = \sigma'_0 \cup
\sigma''_0$.  The set
$$\{f(q) + J^2 M_1 \mid q \in \sigma''_1\} \cup \{f(q) + J^2 M_1
\mid q \text{\ is a path of length\ } 1 \text{\ in \ } P_1\}$$ generates
$J M_1 / J^2 M_1$, and the first of the two subsets is linearly
independent by hypothesis.  Therefore, we may choose a set
$\sigma_1'$ of paths $q$ of length $1$ in $P_1$ such that the images
$f(q) + J^2 M_1$ with $q \in \sigma_1' \cup \sigma_1''$ constitute a
basis for $J M_1/ J^2 M_1$.  Set $\sigma_1 =
\sigma_1'
\cup \sigma_1''$.  Next we find that the union of
$\{f(q) + J^3 M_1 \mid q \in \sigma''_2\}$ with the set of those residue
classes
$f(q) + J^3 M_1$, which correspond to the paths $q$ of length
$2$ in $P_1$ that contain some path in $\sigma'_1$ as a right subpath,
generates $J^2 M_1 / J^3 M_1$; again, the first of the listed sets is
linearly independent by hypothesis.   This permits us to choose a subset
$\sigma'_2$ of the set
$$\{q \mid q \text{\ is a path of length \ } 2 \text{\ in\ } P_1,\  q =
\alpha p \text{\ for an arrow \ } \alpha \text{\ and\ } p
\in \sigma'_1\},$$ so as to obtain a basis $\{f(q) + J^3M_1
\mid q
\in \sigma'_2 \cup
\sigma''_2\}$ for $J^2 M_1 / J^3 M_1$.  We define
$\sigma_2 = \sigma_2' \cup \sigma_2''$, and continue inductively.  Our
construction then guarantees that the resulting set $\sigma =\bigcup_{0
\le l \le L} \sigma_l$ is a skeleton of $M_1$ with the required
properties. \qed
\enddemo

As announced, over a truncated path algebra, any skeleton of a module
determines its syzygies.  More precisely, we obtain:

\proclaim{Theorem 2.3. A known fact} \cite{\BHT, Lemma 5.10}
 If $M$ is a nonzero left
$\la$-module, not necessarily finitely generated, and $\sigma$ any
skeleton of $M$ {\rm {(}}in a suitable projective cover of
$M${\rm {)}}, then
$$\Omega^1(M) \cong \bigoplus_{q\ \S\text{-}\text{critical}} \la q.$$
In particular, $\Omega^1(M)$ is isomorphic to a direct sum of cyclic
left ideals generated by nonzero paths of positive length in
$\la$.  \qed
\endproclaim

The following observation was already used in \cite{\DHL}.

\proclaim{Observation 2.4}  Given a path $q$ of positive length in
$\la$, the cyclic left ideal $\la q$ has finite projective dimension if
and only if the endpoint of $q$ is not precyclic.  As a consequence, a
simple module $\la e / Je$ has finite projective dimension precisely
when it is not precyclic. \qed
\endproclaim

A major asset of truncated path algebras lies in the ease with which
computations implicit in the theory can be carried out graphically via
the {\it layered and labeled graphs\/} of modules as described in
\cite{\dom} and \cite{\menace}.  The following illustration of module
graphs (over a nontruncated finite dimensional algebra $\latilde$) is to
provide an informal reminder of all that is relevant for the present
article.  The algebra we use for this purpose will resurface in Example
10.1.

\definition{Example 2.5}  Let $A = K\Qtilde /
\Itilde$, where $\Qtilde$ is the quiver

$$\xymatrixcolsep{4pc}
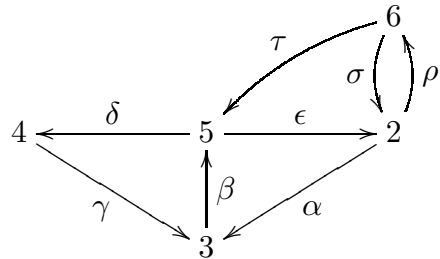
\xymatrix{
 &&6 \ar@/_/[dl]_{\tau} \ar@/_/[d]_{\sigma} \\ 4 \ar[dr]_{\gamma} &5
\ar[l]_{\delta} \ar[r]^{\epsilon} &2
\ar@/_/[u]_{\rho} \ar[dl]^{\alpha} \\
 &3 \ar[u]_{\beta} }$$

\noindent and $\Itilde \subseteq K \Qtilde$ is the ideal generated by
the following relations:
$\alpha \epsilon \beta \alpha$, $\tau \rho$, $\rho \epsilon \beta$,
$\delta \beta \alpha$, $\rho \sigma \rho$,
$\beta \alpha \epsilon \beta$, $\delta \beta \gamma \delta$, $\delta
\tau$, $\epsilon \beta \alpha - \sigma \rho$, and
$\alpha \epsilon - \gamma \delta$.  The following are examples of
layered and labeled graphs of certain left $A$-modules:

$$\xymatrixcolsep{1.0pc}\xymatrixrowsep{2.0pc}
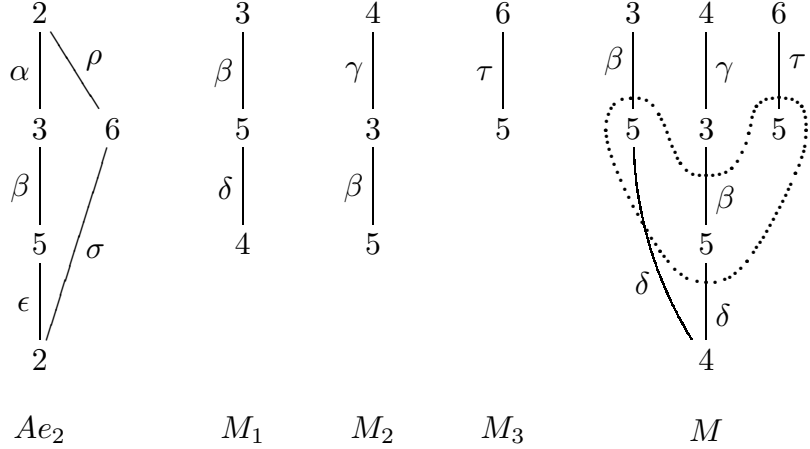
\xymatrix{
2 \edge[d]_{\alpha} \edge[dr]^{\rho} &&&3 \edge[d]_{\beta} &&4
\edge[d]_{\gamma} &&6 \edge[d]_{\tau} &&3 \edge[d]_(0.4){\beta} &4
\edge[d]^{\gamma} &6 \edge[d]^(0.4){\tau} \\
3 \edge[d]_{\beta} &6 \edge[ddl]^{\sigma} &&5 \edge[d]_{\delta}
&&3 \edge[d]_{\beta} &&5 &&5 \edge@/_/[ddr]_(0.6){\delta}
{\save[0,0]+(-3,0);[0,0]+(-3,0) **\crv{~*=<2.5pt>{.} [0,0]+(-3,3)
&[0,0]+(0,4) &[0,0]+(3,3) &[0,1]+(-3,-5) &[0,1]+(0,-5) &[0,1]+(3,-5)
&[0,2]+(-3,0) &[0,2]+(-3,3) &[0,2]+(0,4) &[0,2]+(3,3) &[0,2]+(3,0)
&[0,2]+(3,-3) &[1,1]+(4,-4) &[1,1]+(0,-4) &[1,1]+(-4,-4) &[0,0]+(-3,-3)
}\restore} &3 \edge[d]^(0.6){\beta} &5 \\
5 \edge[d]_{\epsilon} &&&4 &&5 &&&&&5 \edge[d]^(0.6){\delta} \\
2 {\save+<0ex,-5ex> \drop{Ae_2} \restore} &&& {\save+<0ex,-5ex>
\drop{M_1} \restore} && {\save+<0ex,-5ex>
\drop{M_2} \restore} && {\save+<0ex,-5ex> \drop{M_3} \restore} &&&4
{\save+<0ex,-5ex> \drop{M} \restore}
}$$

\noindent The leftmost graph represents the indecomposable projective
left $A$-module $Ae_2$.  Its shape shows that
$\epsilon \beta \alpha - k\sigma \rho$ for some $k \in K^*$; in the
present situation, the scalar is $k = 1$, due to the relations of $A$.
The labeling of the edges is redundant in this example, since the quiver
$\Qtilde$ does not contain multiple edges between any pair of
vertices. (For graphs of the remaining indecomposable projective left
$A$-modules, we refer to Example 9.1.)

The three tree graphs in the center depict uniserial modules $M_i$, $i =
1,2,3$.  For instance, $M_1 \cong A e_3/A\delta \beta$.  If $x_i$ is
a top element of
$M_i$  (see {\it Terminology\/} in Section 1), the rightmost graph
represents the isomorphism class of the module
$$M = \bigl(M_1 \oplus M_2 \oplus M_3 \bigr)/ A \bigl(k_1 \beta x_1 +
k_2 \beta \gamma x_2 + k_3 \tau x_3 \bigr)$$
with $k_1, k_2, k_3 \in
K^*$ (the choice of these scalars clearly does not impinge on the
isomorphism type of $M$); the dotted line enclosing the vertices
representing $\beta x_1, \beta \gamma x_2, \tau x_3$ signifies that
these elements are $K$-linearly dependent, while any two of them are
$K$-linearly independent.  The relation $\delta \beta x_1 \in K^*\delta
\beta \gamma x_2$ in $M$, shown in the graph, is a consequence of the
relation $\delta \tau = 0$.

Clearly, a
graph of a module $X$ need not determine $X$ up to isomorphism, unless
it is a tree.  Conversely, the isomorphism class of
$X$ will typically not lead to a unique graph representing it, but will
do so only once a sequence of top elements of $X$ has been
specified.      An alternate graph of the module $M$ above, for
example, is a disjoint union of two nontrivial subgraphs, thus
displaying decomposability of $M$ at first glance.
\enddefinition

We remark, moreover, that any skeleton $\S$ of a finitely generated
module
$M$ can be pinned down in a visually suggestive format by a finite
forest of tree graphs, one tree for each element in the chosen
sequence of top elements $z_r$ in the distinguished projective cover
$P$ of $M$; we give an illustration below.  The skeleton $\S$ can be
retrieved from its graph as the set of all edge pathis of length $\ge 0$
that start in a vertex representing an element
$z_r \in \S_0$.  Any such graph of a skeleton displays a composition
series of
$M$, recording, from the top down, the simple composition factors in
each of the radical layers
$\SS(M) = \bigl(J^l M / J^{l+1}M \bigr)_{l \ge 0}$. Below, we give
examples of skeleta for two of the modules displayed in Example 2.5.
The left-hand graph depicts the (unique) skeleton of
$Ae_2$, and on the right we exhibit one of the two skeleta of the module
$M$ of that example.  The graphs of the three uniserial modules
$M_i$ coincide with graphs of their skeleta.

$$\xymatrixcolsep{1.0pc}\xymatrixrowsep{2.0pc}
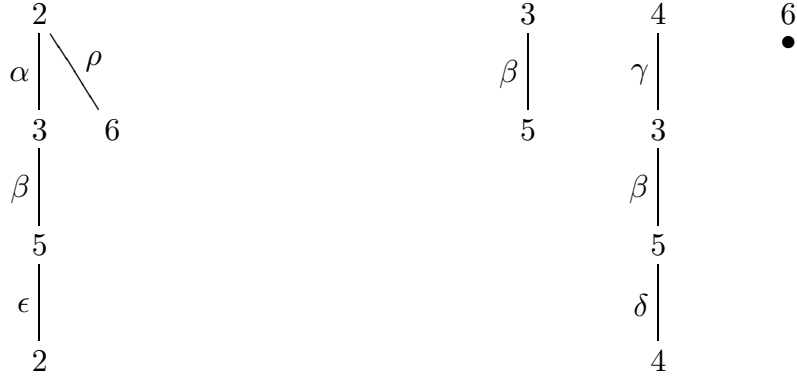
\xymatrix{
2 \edge[d]_{\alpha} \edge[dr]^{\rho} &&&&& &&&3 \edge[d]_{\beta}
&&4 \edge[d]_{\gamma} &&6 \drbl \\
3 \edge[d]_{\beta} &6 &&&& &&&5 &&3 \edge[d]_{\beta} \\
5 \edge[d]_{\epsilon} &&&&& &&& &&5 \edge[d]_{\delta} \\
2 &&& &&&&& &&4
}$$

\head 3. Structure of the $\la$-modules of finite projective dimension.
First installment
\endhead

We continue to let $\la$ denote a truncated path algebra.  The pivotal
result of this section is Theorem 3.1, which exhibits a strong
finiteness property of the categories
$\pinflamod$ and $\Pinflamod$.  This property is responsible for the
first set of structure results which we assemble in Corollary 3.3.
Moreover, Theorem 3.1 will readily yield contravariant finiteness of
$\pinflamod$ and
$\Pinflamod$ in $\lamod$ and $\Lamod$, respectively (Section 4).  To
describe this finiteness property, we  let
$$\eps \ \ \text{be the sum of the idempotents corresponding to the non-precyclic vertices of}\ Q$$
and
observe that, for any $\la$-module $M$, the subspace
$\eps M$ is actually a
$\la$-submodule.  The $\la$-module $\eps M$ always has finite projective
dimension as all of its composition factors are non-precyclic.  The next
theorem will show that the functor
$$F = \la/ \eps J\, \otimes_\la \, - :\, \Lamod \rightarrow
\Lamod,
\ \ M \mapsto M / \eps  J M$$
takes the category $\P^{<
\infty}(\Lamod)$ to a category of finite representation type, namely to
the full subcategory of $\Lamod$, whose objects are direct sums of
copies of the local modules
$\A_i = \la e_i/\eps J e_i$, for $1 \le i \le n$.  Clearly, each
$\A_i$ has finite projective dimension, and equality $\A_i = S_i$ holds
precisely when the corresponding vertex
$e_i$ is non-precyclic.  Moreover, each
$\A_i$ has a tree graph (relative to the top element $e_i +
\eps J e_i$) and is thus determined up to isomorphism by this graph.

Our proof of the next theorem rests on the prerequisites we have
established in the previous section.

\proclaim{Theorem 3.1}  Let $M$ be any
$\la$-module, not necessarily finitely generated.  Then the following
conditions are equivalent:

{\rm (i)} $\pdim M < \infty$.

{\rm (ii)} $F(M)$ is a direct sum of copies of the local modules
$\A_1, \dots, \A_n$.

{\rm (iii)} Given a projective cover $f: P \rightarrow M$, all simple
composition factors of $\Ker(f)$ are

\qquad non-precyclic.

\endproclaim

\definition{Remark} Thus, each object $M$ of $\Pinflamod$ is an
extension of the $\eps \la \eps$-module $\eps M$ by a direct sum of
copies of the (indecomposable projective) $(1-\eps) \la (1 -
\eps)$-modules $\A_i$ corresponding to the precyclic vertices $e_i$,
where all the mentioned modules have a canonical left $\la$-structure,
since
$\eps \la \eps = \la \eps$.  Conversely, every such extension belongs to
$\Pinflamod$. In short, if we arrange the simple modules so that $S_1,
\dots, S_m$ are precyclic and
$S_{m+1},
\dots, S_n$ non-precyclic, we obtain $\eps = e_{m+1} +
\cdots + e_n$ and find
$$\Pinflamod =  \Ext^1_\la \biggl( \Add(\A_1, \dots, \A_m)\, ,\,
\eps \la \eps\text{-}\Mod \biggr),$$  and
$$\pinflamod  = \Ext^1_\la \biggl(\add(\A_1, \dots, \A_m)\, ,\,
\eps
\la \eps\text{-}\mod \biggr).$$
\enddefinition

\demo{Proof of Theorem 3.1} To verify ``(i)$\implies$(ii)", suppose that
$M$ is a nonzero module of finite projective dimension, and let
$f: P \rightarrow M$ be a projective cover.  Then also $\pdim (M /\eps
JM) < \infty$, and
$f$ induces an epimorphism $\overline{f}: P/\eps JP \rightarrow M /\eps
JM$ with kernel contained in $JP/\eps JP$. If
$\overline{f}$ were not an isomorphism, Lemma 2.2 would yield a skeleton
$\sigma$ of $M /
\eps JM$ with a $\sigma$-critical path $q$ ending in a precyclic
vertex:  Indeed, any skeleton $\sigma$ of $M / \eps JM$ would then be
properly contained in a skeleton
$\sigma^+$ of $P/\eps JP$ in $P$.  Any path $q$ of minimal length in
$\sigma^+ \setminus
\sigma$ would be $\S$-critical, and, since all paths in the latter set
difference have positive length and
$\Ker(\overline{f})
\subseteq JP/ \eps JP$ has only precyclic composition factors,
$q$ would be as required.

By Theorem 2.3, this would force a direct summand isomorphic to
$\la q$ into the syzygy of $M / \eps JM$, which, in light of Observation
2.4, would contradict finiteness of the projective dimension of $M /
\eps JM$. Thus
$\overline{f}$ is an isomorphism
$P /\eps JP \cong M / \eps JM$.  Since the former quotient is a direct
sum of copies of the $\A_i$, so is the latter.

To prove ``(ii)$\implies$(iii)'', suppose that $F(M) \cong
\bigoplus_{1 \le i \le n} \A_i^{r_i}$ with $r_i \ge 0$, and let
$f: P \rightarrow M$ be a projective cover.  Then $M/JM \cong F(M)/ J
F(M) \cong \bigoplus_{1 \le i \le n} S_i^{r_i}$, which shows that
$P$ is isomorphic to $\bigoplus_{1 \le i \le n} (\la e_i)^{r_i}$.  Thus
$$P / \eps JP \cong \bigoplus_{1 \le i \le n} \A_i^{r_i} \cong M/
\eps JM$$
by condition (ii).  It follows that $\Ker(f) \subseteq \eps JP$ as
claimed.

The implication ``(iii)$\implies$(i)'' is straightforward.
\qed
\enddemo

Theorem 3.1 moreover shows the category $\P^{< \infty}(\Lamod)$ to have
an unusual closure property under certain types of subobjects.  We call
a submodule $U$ of a module $M$ {\it top-stably embedded in
$M$\/}, in case $U \cap JM = JU$.

\proclaim{Corollary 3.2}  The category $\Pinflamod$ is closed under
top-stably embedded submodules.
\endproclaim

\demo{Proof}  Suppose $M$ in
$\Lamod$ has finite projective dimension and $U \subseteq M$ is a
top-stably embedded submodule. Let
$f_1 : P_1
\rightarrow U$ be a projective cover of $U$.
 By top-embeddedness, we can extend $f_1$ to a projective cover
$f: P = P_1 \oplus P_2 \rightarrow M$.  The argument for Theorem 3.1 now
guarantees that the induced map $\overline{f}: P/ \eps JP
\rightarrow M/ \eps JM$ is an isomorphism, whence
$\overline{f_1}$ is an isomorphism $P_1 /\eps JP_1 \rightarrow U /\eps
JU$.  Consequently, $U / \eps JU$ is a direct sum of copies of the
$\A_i$, which shows
$\pdim U < \infty$.  \qed
\enddemo

For any set $\Psi$ of finitely generated left $\la$-modules, we denote by
$\filt(\Psi)$ the full subcategory of
$\lamod$ having as objects those modules $X$ that have a finite
filtration $X_0 = 0 \subseteq X_1 \subseteq \cdots
\subseteq X_r = X$ with consecutive factors $X_i/X_{i-1} \in
\Psi$.  By $\Filt(\Psi)$ we denote the analogous subcategory obtained by
waiving the requirement that the filtrations with factors in $\Psi$ be
finite, replacing the natural number
$r$ by any ordinal number.

Again let $e_1, \dots, e_m$ be the precyclic vertices of the quiver
$Q$ of $\la$, and $e_{m+1}, \dots, e_n$ the non-precyclic ones; that is,
the idempotent $\eps$ that gives rise to the functor
$F$ above equals $\sum_{m+1
\le i \le n} e_i$.  We record the information provided by Theorem 3.1 in
slightly different form for future reference.

\proclaim{Corollary 3.3.  Structural information on the categories
$\pinflamod$ and\linebreak $\Pinflamod$}
\smallskip

\noindent{\bf A. Simple objects and composition series.}   The set of
simple objects in $\Pinflamod$ is
$$\Psi = \{\A_1, \dots, \A_n\} = \{\A_1, \dots, \A_m, S_{m+1},
\dots, S_n\}.$$  Moreover,
$$\pinflamod = \filt(\Psi) \ \ \text{and} \ \ \Pinflamod = \Filt
(\Psi).$$
Given $M\in \Pinflamod$, the cardinal multiplicities of the simple
$\Pinflamod$-com\-po\-si\-tion factors of
$M$ {\rm(}with respect to an ordinal-indexed composition series{\rm)}
are isomorphism invariants of
$M$.
\smallskip

\noindent{\bf B. Separation of precyclic and non-precyclic portions in
the modules of finite projective dimension.}   Any object $M \in
\Pinflamod$ has a unique largest submodule $U(M) = \eps M$ with the
property that all simple composition factors of $U(M)$ are among
$S_{m+1}, \dots, S_n$, and
$M/U(M)$ is a direct sum of copies of the remaining simple objects in
$\Pinflamod$, namely $\A_1, \dots, \A_m$.  In particular, the simple
composition factors  of $U(M)$ in $\Pinflamod$ and $\Lamod$ coincide, and hence
so do the composition lengths in the two categories, provided that $M$
is finitely generated.
\smallskip

\noindent{\bf C. Separation in the indecomposable projective objects
$\la e_i$.}   For
$i \le m$, we have $\la e_i/U(\la e_i) \cong \A_i$, and given a
non-precyclic simple $S_j$, its multiplicity as a composition factor of
$U(\la e_i)$ is equal to the number of paths of lengths $\le L$ from
$e_i$ to $e_j$.  For $i \ge m+1$, we have
$U(\la e_i) = \la e_i$.
\smallskip

\noindent{\bf D. Finitistic dimensions.} $\lfindim \la =
\Lfindim \la = \max\{\pdim \A_i, \ \pdim S_j \mid i \le m, \ j
\ge m+1\}$.
\endproclaim

\demo{Proof} Let $M$ be any object in $\Pinflamod$.  As was pointed out
above, $U(M) = \eps M$ is then a $\la$-submodule of
$M$, which clearly has a filtration with consecutive factors among
$S_{m+1},
\dots, S_n$.  By Theroem 3.1, $M/U(M)$ is a direct sum of copies of
$\A_1, \dots, \A_m$, and thus, evidently, $M/U(M)$ has a filtration with
factors among the $\A_i$, $i \le m$.   The remaining assertions
 are easy consequences. \qed \enddemo

Corollary 3.3 will be  supplemented in Corollary 5.4, where the
indecomposable (relative) injective objects of
$\Pinflamod$ will be identified.  The combined information will prove
helpful in Section 9, towards exploiting a duality relating the category
$\pinf(\mod$-$\latilde)$ of modules of finite projective dimension over
the strongly tilted algebra
$\latilde$ to the more directly accessible category
$\pinflamod$.

\smallskip
\definition{Remark 3.4.  Standard stratification of truncated path
algebras}  Filtration categories similar to those above have been
studied extensively in the context of quasi-hereditary and standardly
stratified algebras (see, e.g., \cite{\Rin, \Dla, \AgHaLuUn}).  In fact,
the set $\Psi = \{\A_1, \dots, \A_n\} = \{\A_1, \ldots, \A_m, S_{m+1},
\ldots, S_n\}$ in the preceding corollary can be viewed as the set of
{\it standard}
$\Lambda$-modules relative to a pre-order on the set of isomorphism
classes of simple $\Lambda$-modules.  This pre-order is defined by
specifying $S_i \preceq S_j$ if either $e_i$ is precyclic or there exists
a path from $e_i$ to $e_j$ in $Q$, and we write $S_i \prec S_j$ if $S_i
\preceq S_j$ and $S_j \not \preceq S_i$.   It is then easy to see that
$\A_i$ coincides with the standard module $\Delta_i$ -- defined to be
the unique highest-dimensional quotient of $P_i$ having only composition
factors $S_j \preceq S_i$.  Under this preorder on the simples, the
algebra $\la$ is standardly stratified in the sense of
Cline-Parshall-Scott \cite{\CPS}:  Indeed, each projective
$\la e_i$ has a filtration with top factor isomorphic to $\Delta_i$ and
remaining factors isomorphic to
$\Delta_j$ for $S_i \prec S_j$; moreover, the kernel of the canonical
epimorphism $\Delta_i \rightarrow S_i$ has only simple composition
factors $S_j$ with $S_j \preceq S_i$.
The C-P-S standardly stratified  algebras
have been further studied by Frisk
\cite{\Fri}, and some of our results in Section 6 can also
be obtained as applications of the theory he develops.  We continue this
discussion in Remark 6.2.

\enddefinition

\head 4.  Contravariant finiteness of $\pinflamod$ and
$\Pinflamod$
\endhead

\proclaim{Theorem 4.1}  For every truncated path algebra $\la$, the
category $\P^{< \infty}(\lamod)$ is contravariantly finite in
$\lamod$, and $\A_i$ is a minimal (right) $\P^{<
\infty}$-approximation of $S_i$ for $1 \le i \le n$.  In particular, the
minimal approximations of the simple modules are local, and hence
indecomposable.
\endproclaim

\demo{Proof}  As we pointed out in Section 2, contravariant finiteness
of $\pinflamod$ follows from the existence of $\pinflamod$-approximations
of the $S_i$.  Let $\phi: \A_i \rightarrow S_i$ be the canonical
epimorphism, sending the coset $e_i + \eps J e_i$ to
$e_i + J e_i$.  To see that $\phi$ is a $\P^{<
\infty}(\lamod)$-approximation of $S_i$, let
$M$ be any object in $\P^{< \infty}(\lamod)$ and $f$ a nonzero map in
$\hom_\la(M, S_i)$.  Clearly, $f$ factors through the canonical
epimorphism
$\pi$ from $M$ to $M / \eps JM$; say $f = f' \circ
\pi$.  By Theorem 3.1, the latter factor module is isomorphic to a
direct sum of $\A_j$'s. The map $f'$, being an epimorphism onto
$S_i$, clearly factors through a copy of $\A_i$ in any such
decomposition, which guarantees that
$f$ factors through $\phi$.  Minimality of $\phi$ as a
$\pinflamod$-approximation of $S_i$ follows from the indecomposability
of $\A_i$.
\qed \enddemo

It is easy to describe the minimal $\pinflamod$-approximations of
arbitrary finite dimensional $\la$-modules in terms of their projective
covers.  In fact, the description extends to the infinite dimensional
case, thus providing us with minimal
$\Pinflamod$-approximations of arbitrary objects in $\Lamod$.

\proclaim{Theorem 4.2}  Let $M \in \Lamod$, say $M \cong P/C$, where
$P$ is projective and $C \subseteq JP$.  If we identify $M$ with $P/C$,
the canonical map
$\phi: P/ \eps C \rightarrow P/C$ is a minimal (right)
$\P^{< \infty}(\Lamod)$-approximation of $M$.

In particular, $\Pinflamod$ is contravariantly in $\Lamod$.
\endproclaim

\demo{Proof} We note that
$P/\eps C$ has finite projective dimension.  To see that
$\phi$ is a $\Pinflamod$-approximation of $P/C$, let $f: N
\rightarrow P/C$ be any homomorphism with source $N \in
\Pinflamod$.  We clearly do not lose generality in identifying
$N$ with a quotient $Q/D$, where $Q$ is projective and $D
\subseteq JQ$.  To ascertain that $f$ factors through
$\phi$, we consider the following diagram:

$$\xymatrixrowsep{2pc}\xymatrixcolsep{3pc}
\xymatrix{ P \ar[r]^-{\txt{can}} &P/\varepsilon C \ar[r]^-{\phi} &P/C \\
 &&Q/D \dshdar[ul]^-{f''\;} \ar[u]_{f} \\
 &&Q \ar[u]_-{\txt{can}} \dshdar[uull]^{f'\;}  }$$

\noindent Here the map $f'$ making the larger triangle commute exists
due to the fact that $\phi \circ \can$ is a surjection.  To see that
$f'$ induces a map $f''$ rendering the smaller triangle commutative,
 apply Theorem 3.1 to find $\eps D = D$.  Hence the inclusion
$D \subseteq \Ker(f \circ \can)$ implies $f'(D) \subseteq \eps
\Ker(\phi \circ \can) = \eps C$, and we obtain
$f''$ as desired.

To check minimality of $\phi$, let $g$ be an endomorphism of
$P/\eps C$ with $\phi \circ g = \phi$.  Then $g = \id + g'$, where the
image of $g'$ is contained in $C/ \eps C \subseteq J(P/\eps C)$.   In
particular, $g'$ is nilpotent, and hence $g$ is invertible.
\qed
\enddemo

\noindent {\bf Remark concerning skeleta of minimal approximations.} The
description of the minimal $\pinflamod$-approximations in Theorem 4.2
has a simple interpretation in terms of skeleta (introduced in
Definition 2.1).  Given $M = P/C$ as above and a distinguished sequence
of top elements $(z_r)_{r \in R}$ of $P$, any skeleton $\sigma$ of $M$
is contained in the (unique) skeleton $\sigma'$ of $P$ with respect to
these top elements; here $\sigma'$ consists of {\it all} paths in $P$,
that is, of all nonzero elements $pz_r$ where $p$ is a path of length
$\leq L$ in $Q$.  We enlarge $\sigma$ to obtain a skeleton $\sigma''$ of
the $\pinflamod$-approximation $P/\eps C$ of $M$ by adding paths from
$\sigma'$ as follows: The new paths are simply those paths in $\sigma'
\setminus \sigma$ that end in precyclic vertices; clearly these induce a
basis for $C/\eps C$.
\smallskip

The remark on skeleta makes it straightforward to compute the
minimal $\pinflamod$-approximation of any module $M \in \lamod$ from a
minimal projective presentation of $M$.  We illustrate this in Section
5.

Let $A$ be any finite dimensional algebra for the moment. As we
pointed out in Section 2, contravariant finiteness of
$\pinfamod$ implies that the objects of $\pinfamod$ can be described as
the direct summands of the modules
$N$ with the following property:  $N$ has a finite filtration whose
consecutive factors are among the minimal
$\pinfamod$-approximations of the simple left $A$-modules.  In general,
this description of the objects in $\pinfamod$ cannot be sharpened so as
to allow omitting the step of taking direct summands of suitably filtered
modules. Indeed, the category of modules with filtrations as described
need not be closed under direct summands.  However, for truncated path
algebras, it is.  This
was already recorded in Corollary 3.3.
\bigskip

\noindent {\bf A generalization of Theorem 4.1 via a result of Smal\o}
\smallskip

After establishing contravariant finiteness of $\pinflamod$ for
truncated path algebras $\la$, we noticed that Theorem 4.1 permits a
modest generalization by way of a theorem of Smal\o.  In particular,
this approach yields an alternate proof for contravariant finiteness of
$\pinflamod$ when $\la$ is a truncated path algebra.  In the following
$A$ will be an arbitrary Artin algebra.

\proclaim{Theorem 4.3 [\Sma]}  Consider the triangular matrix ring
$$A = \pmatrix \Delta & 0 \\ M & \Gamma \endpmatrix,$$
where
$\Gamma$ and $\Delta$ are Artin algebras and
${}_{\Gamma}M_{\Delta}$ is a bimodule such that
$\pdim{}_{\Gamma} M < \infty.$  Then $\pinf(A$-$\mod)$ is
contravariantly finite in $A$-$\mod$ if and only if
$\pinf(\Delta$-$\mod)$ and $\pinf(\Gamma$-$\mod)$ are contravariantly
finite subcategories of $\Delta$-$\mod$ and
$\Gamma$-$\mod$, respectively.
\endproclaim

Guided by Section 3, we apply this theorem to the following scenario.

\proclaim{Corollary 4.4} Assume that $\eps \in A$ is an idempotent such
that the following hold:

{\rm (i)} $\eps A = A \eps A$;

{\rm (ii)} $\pdim_{A} (A \eps / J \eps) < \infty$ \rm{(in view of (i),
this is equivalent to $\gldim
\eps A \eps < \infty$)};

{\rm (iii)} $\pinf(A/A\eps A$-$\mod)$ is contravariantly finite in
$A/ A \eps A$-$\mod$.

\noindent Then $\pinfamod$ is contravariantly finite in
$A$-$\mod$.
\endproclaim

\demo{Proof}  Condition (i), being equivalent to
$A \eps = \eps A \eps$, implies that $A$ is isomorphic to the triangular
matrix ring $$\pmatrix
 A/ A\eps A & 0 \\ \eps A (1-\eps) & \eps A
\eps \endpmatrix.$$  Clearly, the condition that the corner ring
$\eps A \eps$ has finite global dimension ensures that the remaining
hypotheses of the theorem are satisfied. \qed \enddemo

To derive Theorem 4.1 from Corollary 4.4, we observe:
When $A = \Lambda$ is a truncated path algebra and $\eps$  is the sum of
the non-precyclic primitive idempotents, the quotient
$\Gamma := \la/\la \eps \la$ is again a truncated path algebra.  Its
quiver is given by the full subquiver of the quiver of $\la$ on the
precyclic vertices.  Consequently, any $\Gamma$-module of finite
projective dimension is projective, i.e., $\pinf(\Gamma \text{-} \mod) =
\add({}_{\Gamma} \Gamma)$.  It follows that $\pinf(\Gamma \text{-}
\mod)$ is contravariantly finite in $\Gamma$-$\mod$, with projective
covers serving as $\pinf(\Gamma \text{-} \mod)$-approximations.  Note,
moreover, that the minimal $\pinflamod$-ap\-prox\-i\-ma\-tion of a
precyclic simple $S_i$ can be identified with the minimal $\pinf(\Gamma
\text{-}
\mod)$-ap\-prox\-i\-ma\-tion of $S_i$, which coincides with its
projective cover
$\Gamma e_i = \la e_i/ \eps \la e_i = \A_i$.  Smal\o's proof of Theorem
4.3 makes use of the observation that $\pinf(A \text{-} \mod)$ can be
identified with $\Ext^1_{A} (\pinf(\Gamma \text{-}\mod), \pinf(\Delta
\text{-}\mod) )$.  Specializing to our setting, we hereby recover the
description of $\pinflamod$ given in the remark after Theorem 3.1:
$$\pinflamod = \Ext^1_{\la} \bigl( \add(\la/ \eps \la), \eps \la \eps
\text{-}\mod \bigr).$$

In fact, Theorem 4.3 generalizes to infinite dimensional
modules over a triangular matrix ring $A$, since $\pinf(A \text{-}
\Mod)$ can be identified with
$\Ext^1_{A} (\pinf(\Gamma \text{-}\Mod), \pinf(\Delta
\text{-}\Mod) )$. Hence, in the case of a
truncated path algebra $\la$, Smal\o's arguments will also yield an
alternative proof of contravariant finiteness of the big category
$\Pinflamod$.  The key observations providing the link are as follows.
For a truncated path algebra
$\la$, the categories $\pinf(\la/\eps\la \text{-}\Mod) = \Add(\la/\eps
\la)$ and $\pinf(\eps \la \eps \text{-} \Mod) = \eps \la \eps
\text{-}\Mod$ are contravariantly finite in the (big) module categories
$\la/\eps \la \text{-}\Mod$ and $\eps \la \eps \text{-}\Mod$,
respectively.

\head 5.  Strong tilting modules over truncated path algebras
\endhead

We start with some remarks that hold for arbitrary Artin algebras
$A$.

In \cite{\AuRe}, Auslander and Reiten showed that any contravariantly
finite resolving subcategory $\C$ of
$A$-$\mod$ which is contained in
$\pinfamod$ gives rise to a basic tilting module which is uniquely
determined up to isomorphism by $\C$.  (Following their terminology, we
call a module
$T$ ``basic" if its endomorphism ring is basic in the usual sense.)  In
the special case where the subcategory $\C$ equals
$\pinfamod$, they call this tilting module {\it strong}.   As noted by
Happel and Unger (see \cite{\HaUn}), existence provided, the strong
tilting module plays a distinguished role in the partially ordered set
of all basic tilting objects in
$A$-$\mod$: it is the unique minimal element, the regular left
$A$-module occupying the opposite end of the spectrum. On the other
hand, we do not know of interesting concrete instances, where the strong
tilting module is completely understood, beyond the situation of finite
global dimension of
$A$; in this extreme case the strong tilting module is just the minimal
injective cogenerator.

The purpose of this section is to explore the strong tilting module $T$
over an arbitrary truncated path algebra; in light of Section 4,
existence is guaranteed.  We will see that constructibility of the
minimal
$\pinf$-approximations in
$\lamod$ in this situation allows us to pin down $T$ in terms of the
quiver $Q$ and the Loewy length
$L$ of $\la$.  From these data one can then (with some mild effort)
compute the corresponding tilted algebra
$\latilde =
\End_\la(T)^{\text{op}}$  --  that is, determine quiver and relations
for $\latilde$.

According to \cite{\AuRe}, the basic strong tilting module corresponding
to a contravariantly finite resolving subcategory
$\C$ of $A$-$\mod$ which is contained in $\pinfamod$, is characterized
by the following property:  It is the direct sum of the indecomposable
Ext-injective objects of
$\C$, one from each isomorphism class.  The following lemma provides a
source of Ext-injectives in any contravariantly finite resolving
subcategory
$\C$ of $A$-$\mod$.  It is well known (\cite{\AuSm}), but we include a
short proof for the convenience of the reader.  Yet, this source does
not yield {\it all
\/} Ext-injectives in $\C$, in general, not even in case $\C$ equals the
category
$\pinfamod$ over a left serial string algebra $A$, as Example 5.2 will
show.

Recall that a subcategory $\C$ of $A$-$\mod$ is called {\it resolving\/} if it contains the finitely generated projective left $A$-modules and is closed under extensions and kernels of epimorphisms.

\proclaim{Lemma 5.1} Let $\C$ be a contravariantly finite resolving
subcategory of
$A$-$\mod$, where $A$ is an arbitrary finite dimensional algebra.  Then
the minimal
$\C$-approximation of any finitely generated injective left
$A$-module is Ext-injective in $\C$.
\endproclaim

\demo{Proof} Suppose $E$ is a finitely generated injective left
$A$-module,  $f: B \rightarrow E$  a minimal $\C$-approximation of
$E$, and $0 \rightarrow B @>{\, g\,}>>  X \rightarrow Y
\rightarrow 0$ a short exact sequence in $\C$.  Then $f$ factors through
$g$ since $E$ is injective, that is, $f = h g$ for a suitable
homomorphism
$h : X  \rightarrow E$.  Moreover, our hypothesis on $f$ implies
$h = f j$ for some $j : X \rightarrow B$, because $X$ belongs to
$\C$.  Thus $f = fjg$, and since $f$ is right minimal,
$jg$ is an isomorphism.  This yields splitness of $g$.
\qed
\enddemo

In the following example, $\pinfamod$ is contravariantly finite in
$A$-mod, but the corresponding basic strong tilting module fails to be a
direct summand of the minimal
$\pinfamod$-approximation of the minimal injective cogenerator of
$A$-$\mod$.

\definition{Example 5.2}  Let $A = KQ/I$, where $Q$ is the quiver
$$\xymatrixcolsep{3.5pc}
\xymatrix{ 1 \ar[r]^{\alpha} &2 \ar[r]^{\beta} &3
\ar@(ur,dr)^{\gamma} }$$

\noindent and $I \subseteq KQ$ is generated by
$\gamma\beta\alpha$ and
$\gamma^2$.  Then $A$ is left serial, and hence $\pinfamod$ is
contravariantly finite by \cite{\BHZ} (it is easy to verify
contravariant finiteness directly in this example).  The minimal
$\pinfamod$-approximation of the minimal cogenerator equals $B = (A
e_1)^3 \oplus A e_2$, and thus provides only two isomorphism classes
of indecomposable $\Ext$-injectives in $\pinfamod$.  In particular,
$B$ fails to be a tilting module. \qed \enddefinition

Next we will see that, by contrast, all Ext-injective objects in
$\pinfamod$ are obtained as in Lemma 5.1, provided that $A =
\la$ is a truncated path algebra.   This fact will lead to the announced
description of the basic strong tilting object in $\lamod$.
\bigskip

For the remainder of the section, we again focus on a truncated path
algebra
$\la = KQ/I$ with vertices $e_1, \dots, e_n$, and let
$\A_i$ be the minimal $\pinflamod$-approximation of the simple left
$\la$-module $S_i$, for $1 \le i \le n$, as described in Section 4.
Moreover, we denote by $E(S_i)$ the injective envelope of
$S_i$ and by $\B_i$ the minimal $\pinflamod$-approximation of
$E(S_i)$.  As will be illustrated in the sequel, not only the
$\A_i$, but also the $\B_i$ can be explicitly determined from $Q$
and $L$ by way of Theorem 4.2.  Consequently, Theorem 5.3
will permit us to construct the basic strong tilting module $T \in
\lamod$ from these data.

\proclaim{Theorem 5.3}  Let $S_1, \dots, S_m$ be the precyclic simple
modules, and $S_{m+1}, \dots, S_n$ the non-precyclic ones.  As
before, denote by $\A_i$ the minimal
$\pinflamod$-approximation of $S_i$, and by $\B_i$ the minimal
$\pinflamod$-approximation of
$E(S_i)$, for $1 \le i \le n$. Then the basic strong tilting module in
$\lamod$ is the direct sum
$$T = \bigoplus_{1 \le i \le m} \A_i\ \ \oplus \bigoplus_{m+1
\le i
\le n} \B_i\ .$$

\smallskip

\noindent $\bullet$ Concerning the first subsum:  The categories
$\add\bigl(\bigoplus_{1 \le i \le m} \A_i\bigr)$ and
$\add\bigl(\bigoplus_{1 \le i \le m} \B_i\bigr)$ coincide; that is,
$$\bigoplus_{1 \le i \le m} \B_i\, \cong \,\bigoplus_{1 \le i
\le m} \A_i^{t_i}\, ,$$  for suitable exponents $t_i \ge 1$.  This
direct sum has only precyclic simple composition factors.
\smallskip

\noindent $\bullet$ Concerning the second subsum,
$\bigoplus_{m+1 \le i
\le n} \B_i$: Suppose $i \ge m+1$. Then $\B_i$ is indecomposable, has
the same top as $E(S_i)$, and has exactly one simple composition factor
isomorphic to $S_i$, namely the socle of
$E(S_i)$.  Moreover, every submodule of
$\B_i$ that is not contained in $J\B_i$ contains this copy of $S_i$.  As
for the other simple composition factors of $\B_i$ in the category
$\pinflamod$:  We have,
$U(\B_i) = \eps \B_i = \eps E(S_i)$, and $\B_i/U(\B_i) \cong
\bigoplus_{j
\le m} \A_j^{k_{ij}}$, where $k_{ij}$ is the multiplicity of
$S_j$ as a direct summand of the top
$E(S_i)/J E(S_i) \cong \B_i /J \B_i$.
\endproclaim

\definition{Crucial notation: Summands of the strong basic tilting
module and primitive idempotents in the corresponding tilted algebra}
In the following, we will write the strong basic tilting module $T
\in \lamod$ in the form
$$T = \bigoplus_{1 \le i \le n} T_i,$$  where $T_i = \A_i$ for each
precyclic vertex $e_i$, and $T_i = \B_i$ if $e_i$ is not precyclic.  In
other words, $T_i$ is the minimal
$\pinflamod$-approximation of $S_i$ if $e_i$ is precyclic, and
$T_i$ is the minimal $\pinflamod$-approximation of the injective envelope
$E(S_i)$ otherwise.   Moreover, for each $i \in \{1 \dots, n\}$, we
denote by $\etilde_i$ the canonical projection relative to this
decomposition, followed by the embedding into
$T$, that is, $\etilde_i: T \rightarrow T_i
\hookrightarrow T$.  This yields primitive idempotents
$\etilde_1, \dots, \etilde_n$ in the tilted algebra $\latilde =
\End_\la(T)^{\text{op}}$ which are in obvious one-to-one correspondence
with the $e_i$.
\enddefinition

\demo{Proof of Theorem 5.3} The description of the basic strong tilting
module $T$ in the first claim will follow from Lemma 5.1 once the
two remaining claims have been established.  Indeed, the second
assertion entails that the
$\Ext$-injective object $\bigoplus_{1 \le i \le n} \B_i$ of
$\pinflamod$ has
$n$ pairwise non-isomorphic indecomposable direct summands, and the
final assertion guarantees that each of them occurs
precisely once in the direct sum displayed in the first claim.
\smallskip

 First let $i \le m$.  Clearly,
$E(S_i)$ has only precyclic simple composition factors in this case,
whence the top of $E(S_i)$ is a direct sum of copies of
$S_1, \dots, S_m$.  Consequently, the projective cover $P$ of
$E(S_i)$ is a direct sum $\bigoplus_{1 \le j \le m} (\la e_j)^{m_{ij}}$
for suitable $m_{ij} \ge 0$.  If $E(S_i)
\cong P/C$, then
$\eps C =
\eps P$.  Now apply Theorem 4.2 to obtain $\B_i = P / \eps C =
\bigoplus_{1 \le j \le m} A_j^{m_{ij}}$ as claimed.  Clearly, each of
$\A_1, \dots, A_m$ arises as a direct summand of some
$\B_j$ for $j \le m$.
\smallskip

Next let $i \ge m+1$, and set $\B = \B_i$.  In proving the  description
of $\B$ given in the last part of the theorem, we observe that the
projective cover of
$E(S_i)$ can be described as follows:  Let $(p_r)_{r \in R}$ be  the
different paths of length $\le L$ which end in the vertex $e_i$ and
are maximal with these two properties; by maximality we mean that the
inequality
$\len(p_r) < L$ occurs only in case $p_r$ starts in a source of
$Q$.  If $e(r)$ is the starting point of of $p_r$, then the projective
cover of $E(S_i)$ is
$P = \bigoplus_{r \in R} \la z_r$ with $\la z_r \cong \la e(r)$.  It
is clearly harmless to identify $E(S_i)$ with a factor
module $P/C$, where $C
\subseteq JP$; thus $E(S_i)$ has a sequence of top elements
$x_r = z_r + C$ normed by $e(r)$, respectively.  We use
Theorem 4.2 once again to find $\B = P/\eps C$.  In particular, we
find that the socle
$S_i$ of $E(S_i)$ is  contained in the socle of
$\B$, and is, in fact, the only non-precyclic simple summand in
$\soc(\B)$.  More precisely, we obtain:  $S_i = K p_r (z_r + \eps C)$
for each
$r \in R$, and $p_s(z_r + \eps C) = 0$ for $s
\ne r$.  Moreover, any top element of
$\B$ has the form $z = z' + z''$ where $z'$ is a nontrivial
$K$-linear combination of the residue classes $z_r + \eps C$, and
$z''$ belongs to $J\B$.  For all but the last of the assertions
concerning $\B$, it suffices to show that
$S_i$ is contained in any cyclic submodule of $\B$ which is generated
by a top element $z$.  To verify this containment, we note that, by
construction,
$p_r z'' = 0$ for all $r$, since either
$\len(p_r) = L$ or else
$p_r$ starts in a source of $Q$.  This yields $K p_r z = S_i$, for
any index $r$ for which $z_r + \eps C$ makes a nontrivial appearance in
$z'$, and thus proves the auxiliary statement.   The final assertion is
an immediate consequence of Theorems 3.1 and 4.2.  \qed
\enddemo

Theorem 5.3 allows for multiple occurrences
of the $\A_j$, for $j \le m$, in the direct sum
$\bigoplus_{1 \le i \le m} \B_i$.  Multiplicities
larger than $1$ are a common occurrence in fact.  For the truncated
path algebra $\la_2$ in Example 5.6 below, for instance, $m = 3$ and
$\bigoplus_{1 \le i \le 3} \B_i = \A_1^3 \oplus \A_2 \oplus \A_3$.

Next we extend the description of the injective objects of $\pinflamod$
to the big category $\Pinflamod$ of not necessarily finite dimensional
modules of finite projective dimension.

\proclaim{Corollary 5.4. Information on the category
$\Pinflamod$, second installment}

\noindent Again, let $T = \bigoplus_{1 \le i \le n} T_i$ be the strong
basic tilting module in $\lamod$, where $T_i = \A_i$ for
$1 \le i \le m$, and $T_i = \B_i$ is a minimal
$\pinflamod$-approximation of $E(S_i)$ for $i \ge m+1$.
Then the full subcategory of injective objects of $\Pinflamod$ is equal
to $\Add(T)$.

The indecomposable injective objects $T_1,
\dots, T_m$ are simple in the category
$\Pinflamod$.  For $i \ge m+1$, the non-precyclic submodule $U(T_i)$
of $T_i$ and the precyclic factor module $T_i/U(T_i)$ of $T_i$ {\rm
{(in the sense of Corollary 3.3 B)}}, both objects of $\Pinflamod$, are
described in the final part of {\rm Theorem 5.3}; in particular,
$T_i$ fails to be simple in the category $\Pinflamod$, unless $S_i$ is
injective.
\endproclaim

\demo{Proof}  Since $T$ is a strong tilting module, the injective
objects of
$\pinflamod$ coincide with the objects of $\add(T)$, and, by Theorem 5.3,
$T$ is a cogenerator for $\pinflamod$.  We deduce that $T$ is even a
cogenerator for the category
$\Pinflamod$:  Indeed, any left $\la$-module $M$ of finite projective
dimension is a directed union of its top-stably embedded finitely
generated submodules $M_r$, $r \in R$, all of which belong to
$\pinflamod$ by Corollary 3.2.  Hence $M$ embeds into a direct limit
$\varinjlim E_r$, where $E_r \in
\add(T)$ contains $M_r$ as a submodule.  But by \cite{\HuSm, Observation
3.1}, the latter direct limit belongs to
$\Add(T)$, because $T$ is $\Sigma$-pure injective. Since $\Add(T)$ is
evidently closed under direct summands, this shows that all injective
objects of
$\Pinflamod$ belong to $\Add(T)$.

For the converse, let $M \in \Pinflamod$ and choose $M_r$ in
$\pinflamod$ for $r \in R$, as in the preceding paragraph.  Then
$\Ext^1_\la\bigl(M, T^{(X)}\bigr)$ is an inverse limit of the spaces
$\Ext^1_\la\bigl(M_r, T^{(X)}\bigr)$ $\cong$
$\Ext^1_\la\bigl(M_r,T\bigr)^{(X)}$, all of which are zero.  This shows
$T^{(X)}$ to be injective in $\Pinflamod$ for any index set $X$. \qed
\enddemo

While the structure of the summands $T_1, \dots, T_m$ of the strong
basic tilting module $T \in \lamod$ is transparent, the structure of the
remaining summands $T_{m+1} , \dots, T_n$ is somewhat harder to
visualize from the formal description. The labeled and layered
graphs of the $T_i$ (in the sense of Section 2) permit us to
understand the structure of the $T_i$ in any given example at a
glance. By means of Theorem 4.2, these graphs can be readily obtained
from graphs of the $E(S_i)$, the latter being obvious. We leave the
easy combinatorial proof of the following remark to the reader.

\proclaim{Remark 5.5}  For each $i \le n$, the indecomposable direct
summand
$T_i$ of the basic strong tilting module has a {\rm {(}}unique, up to
isomorphism{\rm {)}} layered and labeled graph without closed edge
paths; in other words, this graph is a tree.  Conversely, $T_i$ is
uniquely determined, up to isomorphism, by this graph. \endproclaim

Instead of spelling out the easy algorithm for constructing these
graphs, we will illustrate the procedure with two examples.
In particular, we will see: Whenever $e_i$ is a non-precyclic vertex,
the layered graph of
$T_i$ may be visualized as a daddy longlegs.  The body is represented by
the socle $S_i$ of $E(S_i)$, and the legs, usually ramified, are in
one-to-one correspondence with the simple summands in the top of
$E(S_i)$.  In the upcoming example, we do not label the arrows in our
quivers, and accordingly omit labels on the edges of the graphs
representing modules; as the considered quivers have no double arrows,
this omission is harmless.   The second of the two specific situations
exhibited will be revisited in Section 10.

\definition{Examples 5.6} Let
$\la_1$ be the truncated path algebra of Loewy length
$2$ based on the quiver $Q_1$ below:

\medskip
$$\xymatrixcolsep{1.0pc}\xymatrixrowsep{2.0pc}
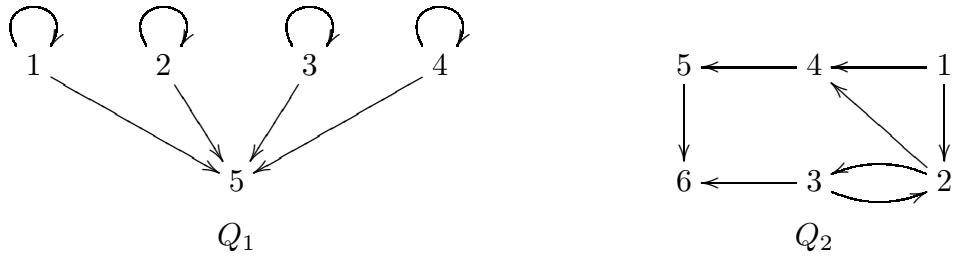
\xymatrix{ 1 \ar@(ul,ur) \ar[drrr] &&2 \ar@(ul,ur) \ar[dr] &&3
\ar@(ul,ur)  \ar[dl] &&4
\ar@(ul,ur) \ar[dlll] &&&&5 \ar[d] &&4 \ar[ll] &&1 \ar[ll]
\ar[d]  \\
 &&&5 {\save+<0ex,-4ex> \drop{Q_1} \restore} &&& &&&&6 &&3
\ar[ll] {\save+<0ex,-4ex> \drop{Q_2} \restore} \ar@/_/[rr] &&2
\ar[ull] \ar@/_/[ll]  }$$

\noindent  The basic strong tilting module in $\la_1$-mod has the
following layered graph:

$$\xymatrixcolsep{1.0pc}\xymatrixrowsep{2pc}
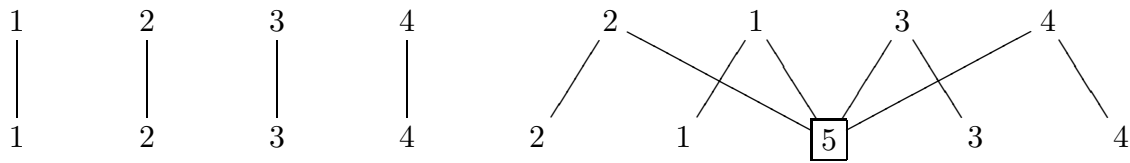
\xymatrix{ 1 \edge[d] &&2 \edge[d] &&3 \edge[d] &&4 \edge[d] &&&2
\edge[dl] \edge[drrr] &&1 \edge[dl] \edge[dr] &&3
\edge[dl] \edge[dr] &&4 \edge[dlll] \edge[dr] \\ 1 &&2 &&3 &&4 &&2 &&1
&&*+[F]{5} &&3 &&4 }$$

\noindent Here the first four trees (from left to right) represent the
summands
$E(S_i) = T_i = \A_i$ corresponding to the precyclic simple modules
$S_1, \dots, S_4$.  The last represents $T_5$, the direct summand
corresponding to the non-precyclic simple module $S_5$; that is, $T_5$
is the minimal
$\pinflamod$-approximation of $E(S_5)$.  In the graph, the socle of
$E(S_5)$  --  the legless body of the spider  -- is highlighted.
\smallskip

Now we consider the truncated path algebra $\la_2$ of Loewy length
$3$ based on the second of the above quivers, $Q_2$. The indecomposable
injective left $\la_2$-modules have graphs:

$$\xymatrixcolsep{0.5pc}\xymatrixrowsep{1.6pc}
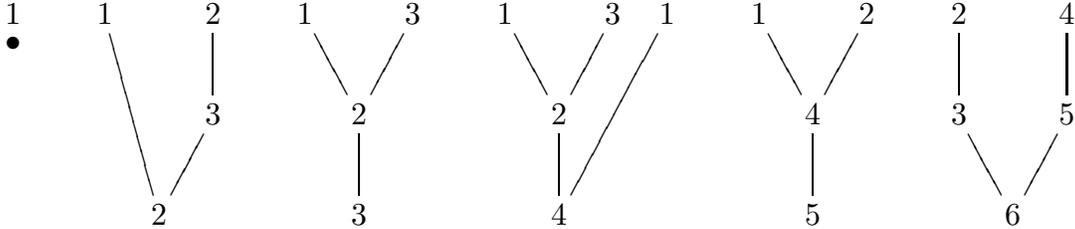
\xymatrix{ 1 \drbl &&1 \edge[ddr] &&2 \edge[d] &&1 \edge[dr] &&3
\edge[dl] &&1 \edge[dr] &&3 \edge[dl] &1 \edge[ddll] &&1
\edge[dr] &&2 \edge[dl] &&2 \edge[d] &&4
\edge[d]  \\
 &&&&3 \edge[dl] &&&2 \edge[d] &&&&2 \edge[d] &&&&&4 \edge[d] &&&3
\edge[dr] &&5 \edge[dl]  \\
 &&&2 &&&&3 &&&&4 &&&&&5 &&&&6 }$$

\noindent Using Theorem 5.3, we obtain the following graph for the basic
strong tilting module in $\la_2$-$\mod$.

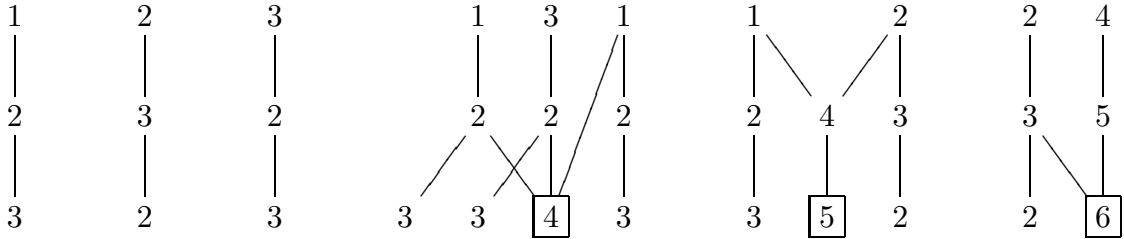
$$\xymatrixcolsep{1.0pc}\xymatrixrowsep{1.6pc}
\xymatrix{ 1 \edge[d] &&2 \edge[d] &&3 \edge[d] &&&1 \edge[d] &3
\edge[d] &1 \edge[ddl]
\edge[d] &&1 \edge[d] \edge[dr] &&2 \edge[dl] \edge[d] &&2
\edge[d] &4
\edge[d]  \\ 2 \edge[d] &&3 \edge[d] &&2 \edge[d] &&&2
\edge[dl] \edge[dr] &2 \edge[dl]
\edge[d] &2 \edge[d] &&2 \edge[d] &4 \edge[d] &3 \edge[d] &&3
\edge[d]
\edge[dr] &5 \edge[d]  \\ 3 &&2 &&3 &&3 &3 &{\save**\frm{-}{4}\restore}
&3 &&3 &{\save**\frm{-}{5}\restore} &2 &&2 &{\save**\frm{-}{6}\restore}
 }$$

\noindent Here the first three trees represent the direct summands
corresponding to the precyclic simples $S_1$, $S_2$,
$S_3$.  The remaining three trees depict the direct summands
corresponding to the non-precyclic simples $S_4$, $S_5$,
$S_6$;  they are the minimal
$\pinflamod$-approximations of $E(S_4)$, $E(S_5)$ , and $E(S_6)$,
respectively, degenerate specimens of spiders in this example.  Again
the bodies of these spiders are highlighted.
\enddefinition

\head 6. Filtrations for the categories $_\la T^\perp$ and
$^\perp (_{\latilde} DT)$
\endhead

In this short section, we develop some structural information on the objects in the
categories of the title.  Such perpendicular categories appear naturally in tilting theory, and are pivotal in transferring information between $\lamod$ and $\latilde$-$\mod$.  In our situation: Since $_\la T$ is strong, we know that $_\la T^\perp = \pinflamod^{\perp}$  (this is essentially the definition of a strong tilting module from \cite{\AuRe}, but see also \cite{\AnTr}).  Our results concerning filtrations of the objects in the targeted subcategories of $\lamod$ and $\latilde$-$\mod$ parallel the
characterization of the objects of $\pinflamod$ in terms of
$\Psi$-filtrations, in Corollary 3.3.

We keep the notational conventions of Section 5, labeling the precyclic
simple left
$\la$-modules $S_1,
\dots, S_m$.   As we saw in Theorem 5.3, the minimal
$\pinflamod$-approximations of the modules in the set
$$\Theta =
\{S_1, \dots, S_m\} \cup \{E(S_{m+1}), \dots, E(S_n)\}$$   are precisely
the indecomposable direct summands $T_i$ of ${}_{\la} T$.  It turns out
that the set $\Theta$ also has a strong impact on the structure of the
objects in the equivalent subcategories $_\la T ^\perp$ of $\lamod$ and
$^\perp(_{\latilde} DT)$ of
$\latilde$-$\mod$.  We write
$\Hom_{\la}(T,\Theta)$ for the set of $\latilde$-modules
$\Hom_{\la}(T,X)$ with $X \in \Theta$, and refer back to Corollary 3.3
for further notation.

\proclaim{Theorem 6.1} We have an equality of subcategories $_\la T
^\perp = \filt(\Theta)$.  Moreover, any $\la$-module $X$ in $_\la
T^\perp$ has a unique largest submodule $V(X)$ with only precyclic
composition factors, and the quotient $X/V(X)$ is a direct sum of copies
of suitable injectives among the $E(S_i)$ for $i \geq m+1$.
 \endproclaim

\demo{Proof} The second statement will follow from our proof of the
first.  Clearly each $E(S_j)$, for $m+1 \leq j \leq n$, belongs to
${}_\la T^{\perp}$.  Next, we observe that $M \in {}_\la T^{\perp}$
whenever $\eps M = 0$.  To see this, we apply the functor $\la/\eps
\la \otimes_{\la} - $ to an extension $0 \rightarrow M
\longrightarrow X \longrightarrow T \rightarrow 0$.  This functor is
exact since $\la/\eps \la \cong (1-\eps)\la_{\la}$ is a projective right
$\la$-module.  Since $T/\eps T$ is projective over
$\la/\eps \la$ (see the remarks following Corollary 4.4), we obtain a
splitting $X/\eps X \rightarrow M$, which yields a splitting of the
original extension upon composition with the canonical map $X
\rightarrow X/\eps X$.  In particular,
$S_i \in {}_\la T^{\perp}$ for each $1 \leq i \leq m$.  Since
${}_\la T^{\perp}$ is extension-closed, we have $\filt(\Theta)
\subseteq {}_\la T^{\perp}$.

For the reverse inclusion, suppose $X \in {}_\la T^{\perp}$, and let
$Y = V(X)$ be the largest submodule of $X$ contained in $\filt(S_1,
\ldots, S_m)$.  The above shows that $Y \in {}_\la T^{\perp}$, and thus
$Z = X/Y \in {}_\la T^{\perp}$ as well, since ${}_\la T^{\perp}$ is
coresolving.  Moreover, maximality of $Y$ ensures that $\soc Z$ has no
precyclic composition factors, i.e., $\soc Z = \eps (\soc Z)$.  We claim
that $Z$ must be injective.  Consider the injective envelope of $Z$
 $$0 \rightarrow Z \longrightarrow E(Z) \longrightarrow W
\rightarrow 0.$$   For each $1 \leq i \leq m$, we clearly have
$\Hom_\la(T_i, Z) = \Hom_\la(T_i, E(Z)) = 0$ since such $T_i$ have only
precyclic composition factors.  By hypothesis,
$\Ext_\la^1(T_i, Z) = 0$, and it follows that $\Hom_\la(T_i, W) = 0$.
Since each precyclic simple occurs in the top of some
$T_i$ with $1 \leq i \leq m$, we must have
$\soc W = \eps (\soc W)$.  We now consider the pullback of the above
extension along the inclusion $\soc W \rightarrow W$.  Since
$\soc W \in \pinflamod$ and $Z \in {}_\la T^{\perp} =
\pinflamod^{\perp}$, the pullback sequence splits.  Hence the inclusion $\soc W
\rightarrow W$ factors through $g : E(Z)
\rightarrow W$.  However, $g(\soc E(Z)) = g(\soc Z) = 0$ implies that
$\soc W = 0$, and hence that $W =0$.  This shows that $Z \in
\add(E_{m+1} \oplus \cdots \oplus E_n)$, and therefore $X \in
\filt(\Theta)$.  \qed \enddemo

Theorem 6.1 moreover yields dual filtrations for the modules in the
subcategory ${}^{\perp}(_{\latilde}DT)$ of $\latilde$-$\mod$:  Indeed, due to Miyashita \cite{\Miya, 1.15}, the adjoint pair of functors
$$\bigl( T\otimes_{\latilde} - \, , \, \Hom_\la(T, - )\bigr)$$
 induces inverse
equivalences
$$^\perp(_{\latilde} DT) \quad\longleftrightarrow\quad (_\la T)^\perp.$$
Consequently,
$$ {}^{\perp}(_{\latilde}DT) = \Hom_{\la}(T,{}_{\la} T^{\perp}) =
\filt\bigl(\Hom_{\la}(T,\Theta)\bigr).$$ In particular, the regular left
$\latilde$-module $\latilde = \Hom_\la(T,T)$ has a $\Hom_\la(T,
\Theta)$-filtration.
 For information about the right regular structure of
$\latilde$ (in a restricted situation), we refer to Theorem 8.2.
\smallskip

\definition{Remark 6.2. Standard stratification of truncated path
algebras}  As noted in Remark 3.4, any truncated path algebra $\la$ is
standardly stratified with respect to the pre-order $$S_i \preceq S_j \
\Leftrightarrow \ \ i\  \text{is\ precyclic\ or\ there\ is\ a\ path\
from\ } i \ \text{to\ } j\ \text{in}\ Q.$$ Corollary 3.3 then asserts
that $\pinflamod$ coincides with the category of $\la$-modules which are
filtered by the standard modules.  Hence, the tilting module
corresponding to this filtration category by \cite{\AuRe} -- usually
called the {\it characteristic} tilting module of $\la$ -- coincides
with the strong tilting module ${}_\la T$.  Furthermore, the modules in
$\Theta$ can be viewed as the {\it proper costandard} $\la$-modules
$\overline{\nabla}_i$, where we define $\overline{\nabla}_i$ to be the
maximal submodule of $E(S_i)$ for which  $\overline{\nabla}_i/ \soc
\overline{\nabla}_i$ has no composition factors $S_j \succeq S_i$.
With these observations, Theorem 6.1 becomes a consequence of the
familiar formula $\filt(\overline{\nabla}) = \filt{\Delta}^\perp$, which
was proved first for quasi-hereditary algebras in \cite{\Rin}, and
recently extended to standardly stratified algebras by Frisk
\cite{\Fri}.  Moreover, Theorem 21 of \cite{\Fri} describes the
subcategory $\Hom_\la(T,\Theta)$ as the subcategory of
$\latilde$-modules that are filtered by the {\it proper standard} left
$\latilde$-modules ${}_{\latilde}\overline{\Delta}_i$, defined dually to
the proper costandard modules.  The fact that ${}_{\latilde}\latilde$
has a suitable filtration with factors among these modules, as pointed
out above, then corresponds to $\latilde$ being right standardly
stratified (with respect to the preorder opposite to $\preceq$).
\enddefinition

\head 7. Dualities induced by strong tilting modules and quivers without
precyclic sources
\endhead

The primary purpose of this section is a characterization of the
truncated path algebras $\la$ with the property that the strong
tilting module
$_\la T$ is also a {\it strong\/} tilting module over
$\latilde = \End_\la(T)^{\text{op}}$.  (Our convention that, for
$f, g \in  \End_\la(T)$, the product
$f\circ g$ stands for ``first apply $g$, then $f$" makes $T$ a right
$\latilde$-module.)  As is to be expected, this situation allows
for particularly effective transfer of information between the
categories
$\pinflamod$ and $\pinf(\mod$-$\latilde)$.  We begin with a review of
known categorical connections induced by tilting and deduce an alternate
argument for a general characterization of strongness in a tilting
module.

For an arbitrary tilting module $_A T_B$ with $B =
\End_A(T)^{\text{op}}$, one not only obtains equivalent pairs of
subcategories of $A$-$\mod$ and
$B$-$\mod$, respectively, but also partial dualities
$A$-$\mod$ $\leftrightarrow$ $\mod$-$B$.  In fact, it follows from
\cite{\Miya, Theorem 3.5} that the functors $\Hom_A( - , T)$ and
$\Hom_B( - , T)$ induce inverse dualities
$$\pinfamod \cap\, {}^{\perp}(_A T)
\quad \longleftrightarrow\quad  \pinfmodb \cap\, {}^{\perp}(T_B). \tag\ddag$$
Consequently, if $X \in
\pinfamod \cap\, {}^{\perp}(_AT)$, applying $\Hom_A(-,T)$ to a minimal
projective resolution of $X$ yields an exact
$\add(T_B)$-coresolution of finite length for $\Hom_A(X,T)$.  Thus we
see that the subcategories linked by the above duality coincide with the
subcategories $\fcog({}_A T)$ and
$\fcog(T_B)$, respectively, consisting of the modules with finite
$\add(T)$-coresolutions, that is, admitting exact sequences of the form
$$0\rightarrow X\rightarrow T_1\rightarrow \cdots \rightarrow T_s\rightarrow 0$$
with $T_i \in \add(T)$.
On the other hand, by \cite{\AuRe, Theorem 5.5b}, $_A T$ is strong if and only if $\pinfamod= \fcog({}_AT)$. Hence, strongness of the tilting module $_A T$ amounts to the equality
 $$\pinfamod \cap\, {}^{\perp}(_A T) = \pinfamod;$$
symmetrically, $T_B$ is strong if and only if $\pinfmodb \cap\,
{}^{\perp}(T_B) = \pinfmodb$.  Miya\-shi\-ta's duality
$(\ddag)$ thus specializes to a duality
$$\pinfamod \quad\longleftrightarrow\quad \pinfmodb$$
in case the tilting bimodule $_AT_B$ is strong on both sides.  This
latter fact, a compelling motivation for exploring strongness of $_A
T_B$, can also be found in
\cite{\AuRe, Proposition 6.6}, where the dual for strong cotilting
modules is stated.  The following convenient criterion for strongness
is stated in dual form in
\cite{\AuRe, Proposition 6.5}, where it is attributed to Auslander and
Green \cite{\AuGr}.  We supply a short alternate proof based on
Miyashita's duality.

\proclaim{Proposition 7.1} \cite{\AuGr}  Let $_A T_B$ be a tilting
module.  Then $T_B$ is a strong tilting module in
$\mod$-$B$ if and only if all simple left $A$-modules embed into
$\soc(_A T)$.
\endproclaim

\demo{Proof} For convenience, we set $\Cal{X} = \fcog({}_A T)$ and
$\Cal{Y} = \fcog(T_B)$.  Suppose that $T_B$ is strong and that
$\Hom_A(S,T) = 0$ for a simple $A$-module $S$.  We shall obtain a
contradiction by showing that $S = 0$.  Since $\Cal{X}$ is
contravariantly finite and resolving, we can find an exact sequence
$0 \rightarrow K \longrightarrow X_1 \longrightarrow X_0
\longrightarrow S \rightarrow 0$ with $X_i \in \Cal{X}$.  Applying
$\Hom_A(-,T)$, we obtain an exact sequence $0 \rightarrow
\Hom_A(X_0,T) \longrightarrow \Hom_A(X_1,T) \longrightarrow Y
\rightarrow 0$ in $\mod$-$B$ with $\Hom_A(X_i, T) \in \Cal{Y}$.  Since
$T_B$ is strong, $\Cal{Y} = \pinfmodb$, and thus $Y \in
\Cal{Y}$.  It follows that $Y \cong \Hom_A(K,T)$, and applying
$\Hom_B(-,T)$ gives a short exact sequence $0 \rightarrow K
\longrightarrow X_1 \longrightarrow X_0 \rightarrow 0$, implying that
$S=0$.

Conversely, suppose that $\Hom_A(S,T) \neq 0$ for all simple
$A$-modules $S$.  Thus there exist nonzero maps from any nonzero
$A$-module $X$ to $T$.  Now suppose $\pdim_B Y < \infty$.  Since
$\Cal{Y}$ contains the projective $B$-modules, we may assume that
$\Omega Y \in \Cal{Y}$.  Hence we have an exact sequence $0
\rightarrow \Hom_A(X,T) \longrightarrow \Hom_A(T_0,T)
\longrightarrow Y \rightarrow 0$ for some $X \in \Cal{X}$, $T_0
\in \add({}_A T)$ and
$f : T_0 \rightarrow X$.  Left-exactness of $\Hom_A(-,T)$ now implies
that $\Hom_A(\Coker(f), T) = 0$, and hence $\Coker(f) = 0$.  Since
$\Cal{X}$ is resolving, $\Ker(f) \in \Cal{X}$, and it follows that $Y
\cong \Hom_A(\Ker(f), T) \in \Cal{Y}$.  Hence $\Cal{Y} =
\pinfmodb$, and $T_B$ is strong.  \qed \enddemo
\bigskip

For the remainder of this section, we return to a truncated path algebra
$\la = KQ/I$ with basic strong tilting module $_\la T$.  In this case,
Proposition 7.1 translates into a straightforward citerion for the
quiver $Q$ equivalent to strongness of $T_{\latilde}$ as a
tilting object in $\mod$-$\latilde$.  By the preceding general
discussion, the equivalence of conditions (1) and (3) in Corollary 7.2
below holds whenever $_AT$ is strong in
$\Amod$, not only in case $A = \la$ is a truncated path algebra.  We
re-emphasize the equivalence in our specialized context for easy
reference in the upcoming applications.

\proclaim{Corollary 7.2}  Let $\la$ be a truncated path algebra with
basic strong tilting module $T \in \lamod$, and $\latilde =
\End_\Lambda(T)^{\text{op}}$.  Then the following statements are
equivalent:
\smallskip

{\rm (1)} $T$ is a strong tilting module in $\mod$-$\latilde$.

{\rm (2)}  The quiver $Q$ of $\la$ has no precyclic source.

{\rm (3)}  The functors
$\Hom_\Lambda(-, T)$ and
$\Hom_{\latilde}(-, T)$ induce dualities between the categories

\quad \ $\pinflamod$ and $\pinf(\mod$-$\latilde)$.
\endproclaim

\demo{Proof}  (1)$\iff$(2):  We refer to the description of $T$ given in
Theorem 5.3.  In light of this theorem, every nonprecyclic simple
$S_i$ occurs in
$\soc_\la (T_i)$, and hence $\bigoplus_{i \ge m+1} S_i$ always embeds
into $\soc_\la (T)$.  On the other hand, a precyclic simple $S_i$ occurs
in $\soc_\la (T)$ precisely when the corresponding vertex
$e_i$ is the endpoint of a path of length
$L$ in $Q$.  That this be satisfied for all precyclic vertices is
tantamount to the requirement that all precyclic vertices be postcyclic,
that is, to non-existence of a precyclic source.

(3)$\iff$(1) was already justified in the more general scenario
considered ahead of Proposition 7.1.
 \qed \enddemo

Corollary 7.2 implies in particular that
$\pinf(\mod$-$\latilde)$ is contravariantly finite in
$\mod$-$\latilde$, whenever the quiver $Q$ of $\la$ is without precyclic
source.  We conjecture that, more strongly,
$\pinf(\mod$-$\latilde)$ is always contravariantly finite in
$\mod$-$\latilde$ for a truncated path algebra $\la$, irrespective of
the shape of the underlying quiver $Q$.

\head 8.  The categories $\pinflatilde$ and $\Pinflatilde$ in the
$\la$-$\latilde$-symmetric situation \endhead

{\it Throughout this section, we assume that the quiver $Q$ has no
precyclic source\/}.   By Corollary 7.2, this places us in the
``$\la$-$\latilde$-symmetric situation" of the section title.

Let $\la$ and $T = \bigoplus_{1 \le i \le n} T_i \in \lamod$ be as
before; in particular, $T_i = \A_i = \la e_i/ \eps J e_i$ for $i \le
m$, and
$T_i$ is the minimal $\pinflamod$-approximation of $E(S_i)$ for
$i \ge m+1$ (cf. Section 5).  Again, $e_1,
\dots, e_m$ are the precyclic vertices in the Gabriel quiver $Q$
of $\la$, and $e_{m+1}, \dots, e_n$ the non-precyclic ones.  Our choice
of corresponding primitive idempotents $\etilde_1, \dots,
\etilde_n$ in the strongly tilted  algebra $\latilde =
\End_\la(T)^{\text{op}}$ was introduced after the statement of Theorem
5.3.  Moreover, we denote by $\Jtilde$ the Jacobson radical of
$\latilde$, and by $\Stilde_i = \etilde_i \latilde/ \etilde_i
\Jtilde$ the simple right $\latilde$-modules corresponding to the
$\etilde_i$.

By our blanket assumption and Corollary 7.2, the functors
$\Hom_\la( -, T)$ and
$\Hom_{\latilde}(-, T)$ induce inverse dualities
$$\pinflamod \leftrightarrow \pinflatilde.$$
Just as in the case of the truncated path
algebra $\la$, the homology of the tilted algebra $\latilde$ is
therefore in turn
governed by a bicoloring of the vertices
$\etilde_i$ of its quiver $\Qtilde$.  This bicoloring of the
$\etilde_i$ is lined up, via $\Hom_\la(T, -)$, with the one that stems
from the placement of the $e_i$ relative to oriented cycles in
$Q$.  By way of caution, we point out that it  does not carry over to
a symmetric placement of the $\etilde_i$ within $\Qtilde$ in general:
Indeed, an idempotent
$\etilde_i$ corresponding to a non-precyclic vertex
$e_i$ may lie on an oriented cycle of $\Qtilde$; see Example 9.1.

Since both $\pinflamod$ and $\pinflatilde$ are contravariantly finite
(in $\lamod$ and $\mod$-$\latilde$, respectively), we know the little
finitistic dimensions on the pertinent sides of
$\la$ and $\latilde$ to coincide with the big finitistic dimensions.
Combining Corollary 7.2 with \cite{\DH}, we thus obtain:

\proclaim{Proposition 8.1} Suppose $\la$ is a truncated path algebra
based on a quiver without precyclic source, and let
$T$ be the corresponding basic strong
$\Lambda$-$\latilde$ tilting bimodule.  Then
$$\Lfindim \la = \lfindim \la = \pdim_\la T = \pdim_{\latilde} T =
\rfindim \latilde = \Rfindim \latilde.  \qed$$  \endproclaim

Theorem 8.5 below will, moreover, permit us to express the right
finitistic dimensions of
$\latilde$ in terms of the simple modules $\Stilde_{m+1}, \dots,
\Stilde_n$ alone.  Namely,
$$\rfindim \latilde = \max\{\pdim \Stilde_j \mid m+1
\le j \le n\}.$$   This follows from the description of composition
series in the category $\pinflatilde$.   This last equality actually
simplifies the situation encountered for the left finitistic dimensions
of
$\la$; indeed, either of the equalities $\lfindim \la =
\max\{\pdim S_j \mid m+1
\le j \le n\}$ and $\lfindim \la = 1 + \max\{\pdim S_j \mid m+1
\le j \le n\}$ can be realized for suitable truncated path algebras; see
\cite{\DHL}.

In view of the duality of Corollary 7.2, the next theorem follows
readily from the mirror symmetric information provided by  Corollaries
3.3 and 5.4.  The reference to precyclic and non-precyclic portions of
modules in
$\Pinflatilde$ in the upcoming results refers to the cycle structure of
the quiver $Q$, not to that of $\Qtilde$.   Paralleling the definition
of the key idempotent $\eps$ in
$\la$, we introduce the idempotent
$$\epstilde = \sum_{m+1 \le i \le n} \etilde_i$$ in $\latilde$.

\proclaim{Theorem 8.2.  The category $\pinflatilde$}  We continue to
assume that the quiver $Q$ has no precyclic source.
\smallskip

\noindent{\bf A. Simple objects and composition series in
$\pinflatilde$.}
 The simple objects of
$\pinf(\mod$-$\latilde)$ are precisely the right
$\latilde$-modules
$$\etilde_i \latilde = \Hom_\la(T_i,T)\ \ \text{for}\ \ 1 \le i
\le m
\ \ \text{and} \ \ \Stilde_i = \Hom_\la(S_i,T)\ \ \text{for}\ \ m+1
\le i \le n.$$ Moreover, $\etilde_i \Jtilde \ne 0$ for $i \le m$.

In particular, a simple right $\latilde$-module
$\Stilde_i$ has finite projective dimension if and only if $i
\ge m+1$.
\smallskip

\noindent{\bf B. Heredity property: Separation of precyclic and
non-precyclic portions in the objects of $\pinflatilde$.} Each module
$\Mtilde$ in $\pinf(\mod$-$\latilde)$ has a unique subobject
$U(\Mtilde)$ maximal with respect to being a direct sum of copies of the
projective modules
$\etilde_i \latilde$ with $i \le m$.  All simple composition factors of
$\Mtilde / U(\Mtilde)$ are among $\Stilde_{m+1}, \dots,
\Stilde_n$.

  In particular, $U(\Mtilde)$ equals $\Mtilde (1 - \epstilde)
\latilde$, the submodule of
$\Mtilde$ generated by all elements $x$ with $x \etilde_i = x$ for some
$i \le m$, and this module is projective.  Moreover,  if $M \in
\pinflamod$ and $\Mtilde \cong \Hom_\la(M, T)$, the composition length
of $\Mtilde/ U(\Mtilde)$ in $\mod$-$\latilde$ coincides with that of the
$\la$-module $U(M)$ of {\rm Corollary 3.3}.

\smallskip

\noindent{\bf C. Separation in the indecomposable projective objects
$\etilde_i \latilde = \Hom_\la(T_i, T)$.}     For $i
\le m$, the submodule
$U(\etilde_i \latilde)$ of Part B equals
$\etilde_i \latilde$.  Now suppose $i \ge m+1$.  In this case,
$U(\etilde_i \latilde) \cong \bigoplus_{j \le m} (\etilde_j
\latilde)^{k_{ij}}$, where $k_{ij}$ is the multiplicity of
$S_j$ in
$T_i / J T_i$; thus, $k_{ij}$ equals the number of those paths of length
$L$ in $Q$, which start in the precyclic vertex $e_j$ and end in $e_i$.
\smallskip

\noindent{\bf D. The indecomposable injective objects of
$\pinflatilde$.}  These are the $\latilde$-modules
$\Etilde_i =
\Hom_\la(\la e_i, T) \cong e_iT_{\latilde}$.  Each $\Etilde_i$ has a
unique simple subobject in the category
$\pinf(\mod$-$\latilde)$.  This subobject is $\etilde_i
\latilde$ in case $i \le m$, and equals $\Stilde_i$ in case $i
\ge m+1$.  In particular,
$U(\Etilde_i) = \etilde_i \latilde$ for $i \le m$,  and
$U(\Etilde_i) = 0$ for $i \ge m+1$.
\endproclaim

\demo{Proof}  In light of the duality of Corollary 7.2, it is routine to
translate the assertions of Corollaries 3.3 and 5.4 into statements
concerning $\pinflatilde$.  Only the claim that $\etilde_i \Jtilde
\ne 0$ for $i \le m$ (under A) requires further backing.  Since $e_i$ is
a precyclic vertex and $Q$ is free of precyclic sources, $e_i$ is also
postcyclic.  In particular, there exists an arrow from
$e_j$ to $e_i$ for a suitable index $j$, possibly equal to $i$. Since
$e_j$ is clearly again precyclic, we have $j \le m$.  This, in turn,
places a composition factor
$S_i$ into $J T_j$, and thus gives rise to a nonzero map in
$\Hom_\la(T_i, T_j)$ which fails to be an isomorphism.  \qed
\enddemo

Some of the mirror-symmetric structure statements for objects in
$\pinflamod$ and $\pinflatilde$ can be pushed beyond the categorical
level as follows:  Namely, if
$M \in \pinflamod$ and $\Mtilde = \Hom_\la(M,T)$, then
$\eps\, M = M$ is equivalent to $\Mtilde\, \epstilde =
\Mtilde$.  More generally, the equality
$\bigl( \Mtilde / U(\Mtilde) \bigr )\, \epstilde = \Mtilde / U(\Mtilde)$
is tantamount to $\eps\, U(M) = U(M)$.  On the other hand, such
non-categorical dual statements are not consistently available:  While
$\eps\, (M/ U(M)) = 0$, we find $ U( \Mtilde)\, \epstilde \ne 0$, in
general.  In particular, while $\eps \A_i = 0$ for $i
\le m$, the corresponding projective $\latilde$-module
$\widetilde{\A_i} \cong \etilde_i \latilde$ typically has composition
factors of finite projective dimension in its radical (see Example 9.1).

This latter fact is contrasted by the following ``non-heredity"
condition for the $\etilde_i \latilde$ with $i
\le m$.

\proclaim{Corollary 8.3}  Whenever $\Stilde_i$ has infinite projective
dimension, that is, when $i \le m$ in our ordering of the vertices, the
projective module $\etilde_i \latilde$ has no proper nonzero submodule
of finite projective dimension. \qed \endproclaim

We now turn our attention to the big category $\Pinflatilde$.
Essentially, all of the structure results above carry over, but require
additional argumentation, as we are leaving the stage of the duality
$\pinflamod \leftrightarrow \pinflatilde$.

\proclaim{Proposition 8.4} Let $\Mtilde \in \Pinflatilde$ be such that
all simple summands of $\Mtilde/\Mtilde \Jtilde$ in
$\Mod$-$\latilde$ have infinite projective dimension.   Then
$\Mtilde$ is projective, that is
$$\Mtilde\cong  \bigoplus_{i \le m} \bigl(\etilde_i
\latilde\bigr)^{(\tau_i)}$$ for suitable cardinal numbers
$\tau_i$.
\endproclaim

\demo{Proof}  We first observe that the hypothesis on $\Mtilde /
\Mtilde \Jtilde$ is tantamount to the equality $\Mtilde =
\Mtilde(1 -
\epstilde) \latilde$.   So, if $\Mtilde$ is finitely generated, then
$\Mtilde = U(\Mtilde)$ is projective and has the required format by the
preceding theorem.

Now we drop the extra hypothesis on $\Mtilde$.  From
\cite{\HuSm, Theorem 4.4}, it follows that $\Mtilde$ is a direct limit
of a directed system of finitely generated
$\latilde$-modules of finite projective dimension; let
$\Mtilde_r$, $r \in R$, be the members of such a system.  Since
$\Mtilde = \Mtilde (1 -
\epstilde) \latilde$, it is harmless to assume that also
$\Mtilde_r = \Mtilde_r (1 - \epstilde) \latilde$ for all $r \in R$.   As
we already saw, this ensures that all of the
$\Mtilde_r$ are projective of the correct format and, flatness being the
same as projectivity over a finite dimensional algebra, we find that
their direct limit, $\Mtilde = \varinjlim
\Mtilde_r$, is projective as well.  It is clear that only  projectives
with tops in $\{\Stilde_1,
\dots, \Stilde_m\}$ will occur as indecomposable direct summands of
$\Mtilde$.  \qed
\enddemo

\proclaim{Theorem 8.5.  The category $\Pinflatilde$}  We continue to
assume that the quiver $Q$ has no precyclic source.
\smallskip

\noindent{\bf A. Composition series.}  Each object $\Mtilde \in
\Pinflatilde$ has an ordinal-indexed composition series with consecutive
factors among the relative simple objects
$\etilde_i \latilde$ for $i \le m$ and $\Stilde_i$ for $i \ge m+1$. The
cardinalities in which these factors occur are isomorphism invariants of
$\Mtilde$.
\smallskip

\noindent{\bf B. Heredity property: Separation of precyclic and
non-precyclic portions in the objects of $\Pinflatilde$.} Every object
$\Mtilde$ in $\Pinflatilde)$ has a unique submodule
$U(\Mtilde)$ which is maximal with respect to being a direct sum of
copies of the projectives $\etilde_i \latilde$ with
$i \le m$, and all simple composition factors of
$\Mtilde / U(\Mtilde)$ are among $\Stilde_{m+1}, \dots,
\Stilde_n$.

In particular, $U(\Mtilde) = \Mtilde (1 - \epstilde)
\latilde$ is generated by all the elements $x \in \Mtilde$ with $x
\etilde_i = x$ for some $i
\le m$, and this module is projective.
\smallskip

\noindent{\bf C. The injective objects of the category
$\Pinflatilde$.} The injective objects of
$\Pinflatilde$ coincide with the objects in $\Add(\bigoplus_{1
\le i
\le n} \Etilde_i) = \Add(T_{\latilde})$.
\endproclaim

\demo{Proof}  We first prove part B.  Let $\Mtilde \in
\Pinflatilde$.  Then $U(\Mtilde) = \Mtilde (1 -
\epstilde) \latilde$ clearly again belongs to $\Pinflatilde$.  By
Proposition 8.4, we infer that $U(\Mtilde)$ is a direct sum of copies of
suitable
$\etilde_i
\latilde$ with $i \le m$.  By construction $\bigl(\Mtilde / U(\Mtilde)
\bigr) \epstilde = \Mtilde / U(\Mtilde)$, whence all simple composition factors
of this quotient are among $\Stilde_{m+1},
\dots,
\Stilde_n$.  The claim regarding composition series in
$\Pinflatilde$ follows.

Part C can be established in analogy with Corollary 5.4, given that any
$\Mtilde \in \Pinflatilde$ is a directed union of submodules
$\Mtilde_r \in \pinflamod$. \qed
\enddemo

\definition{Remark 8.6. Standard stratification of truncated path
algebras}  If $Q$ has precyclic sources, we can still make use of the
tilting duality $(\ddag)$ displayed in Section 7 to obtain a duality
$$ \pinflamod \leftrightarrow \pinf(\mod \text{-}\latilde) \cap
{}^{\perp}(T_{\latilde}).$$ Hence, the conclusions of Theorem 8.2 are
still valid for the category $\pinf(\mod \text{-}\latilde) \cap
{}^{\perp}(T_{\latilde})$.  Although this category is a proper
subcategory of $\pinflatilde$ when $Q$ has a precyclic source, it always
contains the projective right $\latilde$-modules.  In particular, we see
that $\latilde_{\latilde}$ has a filtration with factors belonging to
the set $$\tilde{\Psi} :=  \Hom_\la(\Psi, T) = \{\etilde_1 \latilde,
\ldots, \etilde_m \latilde, \Stilde_{m+1}, \ldots, \Stilde_n\},$$ with
$\Psi$ as in Corollary 3.3.  Since $\Psi$ is the set of standard
$\la$-modules, it follows from \cite{\Fri, Proposition 23} that
$\tilde{\Psi}$ consists precisely of the standard right modules of
$\latilde$ relative to the opposite of the pre-order on the simple
$\la$-modules: $\Stilde_i
\preceq \Stilde_j \Leftrightarrow S_j \preceq S_i$.  That is to say, the
module $\Hom_\la(\A_i,T)$ is the largest quotient of $\etilde_i
\latilde$ all of whose composition factors $\Stilde_j$ satisfy
$\Stilde_j \preceq \Stilde_i$.  As we observed earlier, $\latilde$ is
right standardly stratified, and thus each projective
$\etilde_i\latilde$ has a filtration with top factor isomorphic to
$\Hom_\la(\A_i, T)$ and subsequent factors isomorphic to $\Hom_\la(\A_j,
T)$ for $j$ with $S_j \prec S_i$.
\enddefinition

\head 9. Examples \endhead

We first give an example based on a quiver $Q$ without precyclic
source.

\definition{Example 9.1}  Let $Q$ be the quiver

$$\xymatrixcolsep{1.0pc}\xymatrixrowsep{2.0pc}
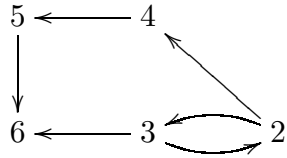
\xymatrix{ 5 \ar[d] &&4 \ar[ll]  \\ 6 &&3 \ar[ll] \ar@/_/[rr] &&2
\ar[ull] \ar@/_/[ll]  }$$

\noindent and $\la$ the truncated path algebra $KQ/I$, where $I$ is
generated by all paths of length $3$.  The strong basic tilting module
$T = \bigoplus_{2 \le i \le 6} T_i$ is determined up to isomorphism by
its graph, namely the following forest:

$$\xymatrixcolsep{0.8pc}\xymatrixrowsep{1.6pc}
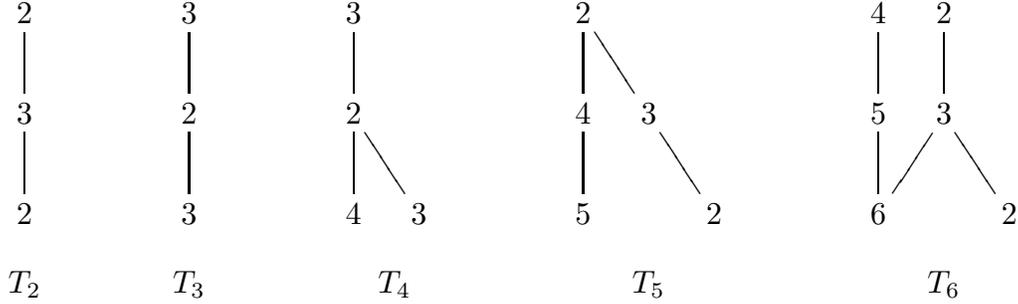
\xymatrix{ 2\edge[d] &&&3 \edge[d] &&&3 \edge[d] &&&&2
\edge[d] \edge[dr] &&&&&4
\edge[d] &2 \edge[d] \\ 3 \edge[d] &&&2 \edge[d] &&&2 \edge[d]
\edge[dr] &&&&4 \edge[d] &3
\edge[dr] &&&&5 \edge[d] &3 \edge[dl] \edge[dr] \\ 2 {\save+<0ex,-5ex>
\drop{T_2} \restore} &&&3 {\save+<0ex,-5ex>
\drop{T_3} \restore} &&&4 {\save+<3ex,-5ex> \drop{T_4}
\restore} &3 &&&5  &{\save+<0ex,-5ex> \drop{T_5} \restore} &2 &&&6 &
{\save+<0ex,-5ex>
\drop{T_6} \restore} &2 }$$

\noindent As in Example 5.6(2), one obtains $T$ as follows:  First one
identifies the cyclebound vertices of
$Q$, namely $e_2$ and $e_3$ in this case, whence $T_i = \A_i$ is the
minimal $\pinflamod$-approximation of $S_i$ for $i = 2,3$, and $T_i$ is
the minimal $\pinflamod$-approximation of the injective envelope
$E(S_i)$ for $i = 4,5,6$.

We will give some detail on the computation of quiver and relations of
$\latilde = \End_\la(T)^{\text{op}}$ in the more challenging Example
9.2, and leave the present case as an exercise.  In the present
example,
$\End_\la(T)$ coincides with the algebra $A$ presented in Example 2.5.
Note that the Gabriel quiver of $A$  --  it is displayed in
2.5   --  contains no multiple edges, whence we may suppress the
labeling of the edges of the graphs of the indecomposable projective
right
$\latilde$-modules (= indecomposable projective left
$\End_\la(T)$-modules) without losing information.  They are as
follows:

$$\xymatrixcolsep{0.8pc}\xymatrixrowsep{1.6pc}
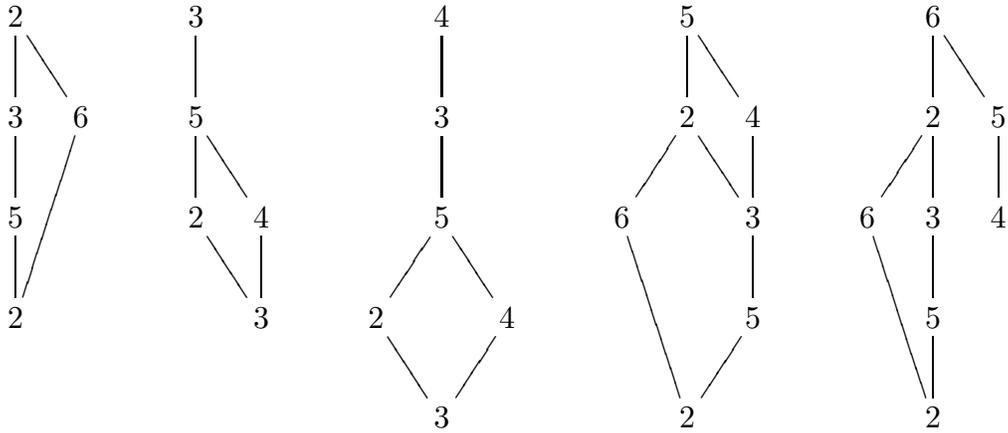
\xymatrix{ 2 \edge[d] \edge[dr] &&&3 \edge[d] &&&&4 \edge[d] &&&&5
\edge[d]
\edge[dr] &&&&6 \edge[d] \edge[dr] \\ 3 \edge[d] &6 \edge[ddl] &&5
\edge[d] \edge[dr] &&&&3 \edge[d] &&&&2 \edge[dl]
\edge[dr] &4 \edge[d] &&&2 \edge[dl] \edge[d] &5
\edge[d] \\ 5 \edge[d] &&&2 \edge[dr] &4 \edge[d] &&&5
\edge[dl] \edge[dr] &&&6
\edge[ddr] &&3 \edge[d] &&6 \edge[ddr] &3 \edge[d] &4 \\ 2 &&&&3 &&2
\edge[dr] &&4 \edge[dl] &&&&5 \edge[dl] &&&5
\edge[d] \\
 &&& &&&&3 &&&&2 &&&&2 }$$

\noindent Since the quiver $Q$ has no precyclic source, the bimodule
$_\la T _{\latilde}$ is strong on both sides (Section 7), and Theorem
8.5 provides us with information on the big category $\Pinflatilde$:
The simple objects are
$\etilde_i \latilde$ for $i =2,3$ and $\Stilde_i = \etilde_i
\latilde/ \etilde_i \Jtilde$ for $i = 4,5,6$.  Moreover, the objects in
$\Pinflatilde$ are precisely the modules $\Mtilde$ with the following
structure:  each has a unique largest submodule
$U(\Mtilde)$ which is a direct sum of copies of
$\etilde_i \latilde$, $i = 2,3$, and $\Mtilde/ U(\Mtilde)$ has only
compostion factors among $\Stilde_4,
\Stilde_5, \Stilde_6$.  In particular, we glean from the above graphs
that $U(\etilde_4 \latilde) \cong \etilde_3
\latilde$ and $\etilde_4 \latilde / U(\etilde_4
\latilde) \cong \Stilde_4$; similarly, $U(\etilde_5 \latilde)
\cong
\etilde_2 \latilde$ and $\etilde_5 \latilde / U(\etilde_5
\latilde)$ is the two-dimensional uniserial module with top
$\Stilde_5$ and socle $\Stilde_4$;  finally, $U(\etilde_6
\latilde) \cong
\etilde_2 \latilde$ with $\etilde_6 \latilde / U(\etilde_6
\latilde)$ is the three-dimensional uniserial module with radical
layering $\bigl(\Stilde_6, \Stilde_5,
\Stilde_4\bigr)$.  The injective objects of $\Pinflatilde$ are the
direct sums of
$\Etilde_i = \Hom_\la(T, T_i)$ for $2 \le i \le 6$.  The general theory
moreover tells us that, for $i = 2,3$, the module $\Etilde_i$ is an
essential extension of
$U(\Etilde_i) = \etilde_i \latilde$, and for $i = 4,5,6$, the only
composition factors of $\Etilde_i$ are among
$\Stilde_j$, $4 \le j \le 6$. $\qed$
\enddefinition

Our final example has a single precyclic source.  In fact, the quiver of
the truncated path algebra we will consider next results from that of
Example 9.1 by the addition of a single source labeled $1$.  As we
will see, this destroys the duality encountered in the previous
example.  In particular, the simple left
$\la$-modules of finite projective dimension are no longer in one-to-one
correspondence with the simple right
$\latilde$-modules of finite projective dimension.  The more intricate
duality theory that governs the general situation will be explored in a
sequel to this paper.

\definition{Example 9.2}  This time, $\la$ is the truncated path
algebra of Loewy length $3$ based on the quiver $Q = Q_2$ of Example
5.6 for which we already computed the basic strong tilting module
$T = \bigoplus_{1 \le i \le 6} T_i$.  Again we exhibit the quiver of
$\latilde = \End_\la(T)^{\text{op}}$, followed by  graphs of the
indecomposable projective right
$\latilde$-modules.

For ease of reading, we display the opposite of
$\Qtilde$, namely, the quiver of $\End_\la(T)$.  Subsequently, we
present the graphs of the indecomposable projective right
$\latilde$-modules (= indecomposable projective left
$\End_\la(T)$-modules).   We label edges only where the two arrows from
vertex $4$ to vertex $1$ might otherwise lead to ambiguities.

$$\xymatrixcolsep{4pc}
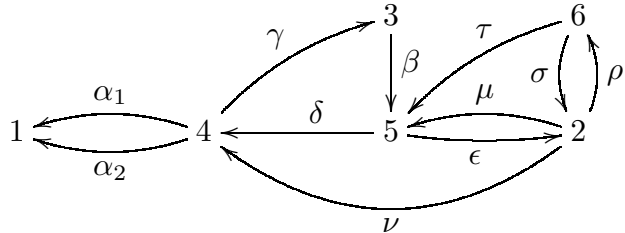
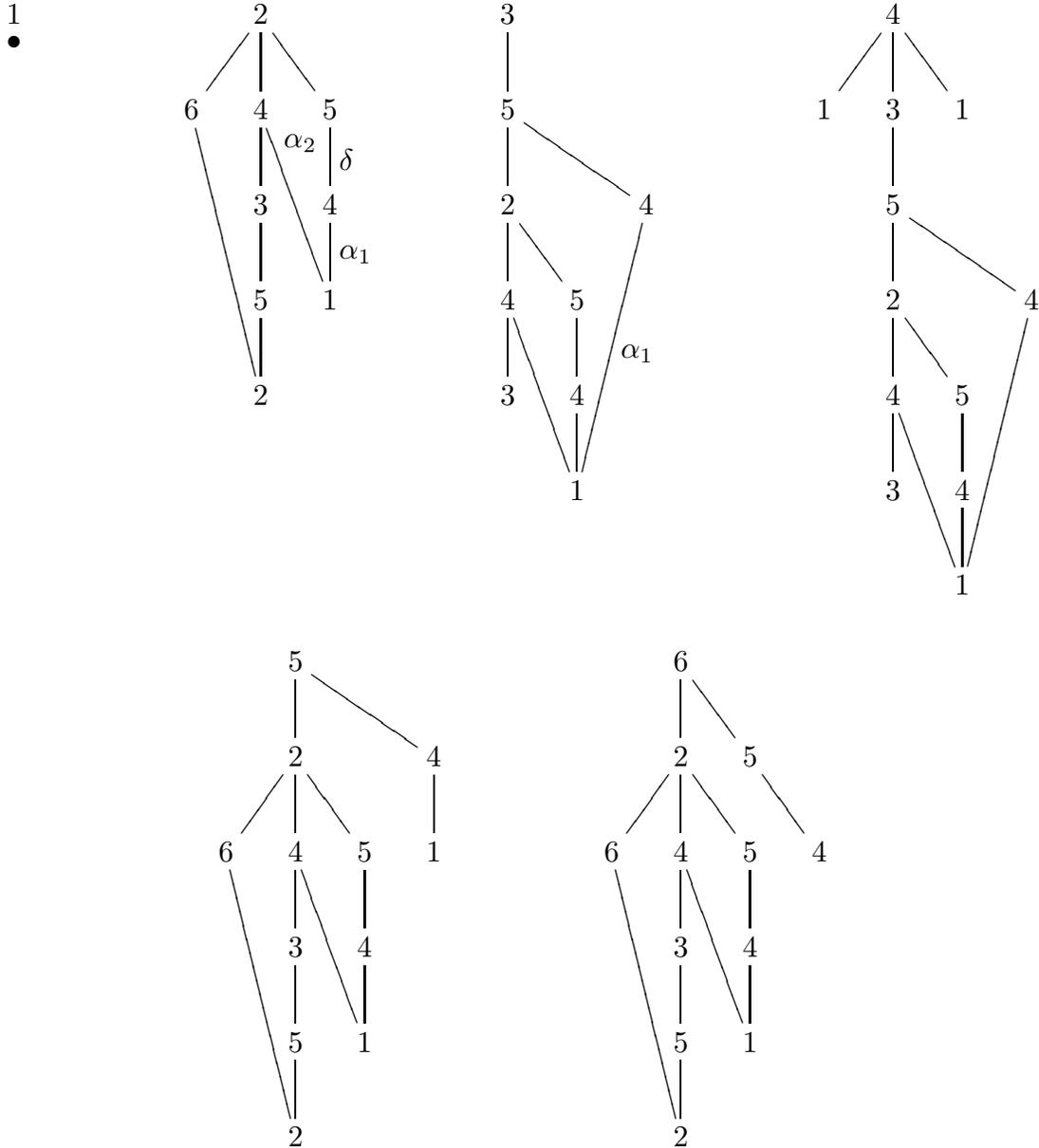
\xymatrix{
 &&3 \ar[d]^(0.4){\beta} &6 \ar@/_/[dl]_(0.4){\tau}
\ar@/_/[d]_{\sigma}  \\  1 &4 \ar@/_/[l]_{\alpha_1} \ar@/^/[l]^{\alpha_2}
\ar@/^/[ur]^{\gamma} &5 \ar[l]_(0.4){\delta}
\ar@/_0.2pc/[r]_(0.45){\epsilon} &2 \ar@/_/[l]_{\mu}
\ar@/_/[u]_{\rho} \ar@/^2pc/[ll]^{\nu} }$$
\bigskip

$$\xymatrixcolsep{1.0pc}\xymatrixrowsep{1.6pc}
\xymatrix{ 1 \drbl &&&&2 \edge[dl] \edge[d] \edge[dr] &&&&3 \edge[d]
&&&&&&4 \edge[dl]
\edge[d] \edge[dr]  \\
 &&&6 \edge[dddr] &4 \edge[d] \edge[ddr]^(0.2){\alpha_2} &5
\edge[d]^{\delta} &&&5
\edge[d] \edge[drr] &&&&&1 &3 \edge[d] &1  \\
 &&&&3 \edge[d] &4 \edge[d]^{\alpha_1} &&&2 \edge[d] \edge[dr] &&4
\edge[dddl]^{\alpha_1} &&&&5
\edge[d] \edge[drr]  \\
 &&&&5 \edge[d] &1 &&&4 \edge[d] \edge[ddr] &5 \edge[d] &&&&&2 \edge[d]
\edge[dr] &&4 \edge[dddl]  \\
 &&&&2 &&&&3 &4 \edge[d] &&&&&4 \edge[d] \edge[ddr] &5 \edge[d]  \\
 &&&&& &&& &1 &&&&&3 &4 \edge[d]  \\
 &&&& &&&& &&&& &&&1 }$$

$$\xymatrixcolsep{1.0pc}\xymatrixrowsep{1.6pc}
\xymatrix{
 &5 \edge[d] \edge[drr] &&&&&&6 \edge[d] \edge[dr]  \\
 &2 \edge[dl] \edge[d] \edge[dr] &&4 \edge[d] &&&&2 \edge[dl] \edge[d]
\edge[dr] &5 \edge[dr]  \\ 6 \edge[dddr] &4 \edge[d] \edge[ddr] &5
\edge[d] &1 &&&6 \edge[dddr] &4
\edge[d] \edge[ddr] &5 \edge[d] &4  \\
 &3 \edge[d] &4 \edge[d] &&&&&3 \edge[d] &4 \edge[d]  \\
 &5 \edge[d] &1 &&&&&5 \edge[d] &1  \\
 &2 &&&&&&2 }$$

\noindent Note that the simple module
$\Stilde_1$ corresponding to the precyclic vertex $e_1$ of $Q$ is
projective, that is, its graph consists only of the vertex
$1$.  Generators for $\Itilde$ such that $\latilde
\cong K\Qtilde/\Itilde$ can be read off the graphs of the indecomposable
projective right modules above, as
$\Itilde$ can be generated by monomial relations and
binomial relations of the form $p - q$, where $p$ and $q$ are paths in
$\Qtilde$.

For our (partial) glossary, we refer the reader back to the graph of the
basic strong tilting module $T = \bigoplus_{1 \le i \le 6} T_i$ which
was displayed in Example 5.6(2).  For $i = 1,2$, let
$\alpha_i: T_4 \rightarrow T_1$ be the epimorphism sending the top
element $x_i$ of $T_4$ shown in the graph below to the top of
$T_1$, and sending the $x_j$ for $j \ne i$ to zero.  Then
$\alpha_1$, $\alpha_2$ are clearly
$K$-linearly independent modulo $\etilde_4 \Jtilde^2
\etilde_1$; thus they yield two arrows $\etilde_4 \rightarrow
\etilde_1$ in the quiver of $\End_\la(T)$.  Next, $\beta \in
\Hom_\la(T_3,T_5)$ is chosen so as to send a top element of
$T_3$ to an element in $e_3 JT_5 \setminus J^2T_5$; any such choice lies
outside $\etilde_3 \Jtilde^2
\etilde_5$.  The map $\gamma
\in \Hom_\la(T_4, T_3)$ is the obvious epimorphism with kernel
$\la x_1 + \la x_2$.  We content ourselves with spelling out a few not
quite so obvious additional choices of arrows.  Namely, we have two
$K$-linearly independent homomorphisms
$T_2 \rightarrow T_4$: one of them is $\nu$, sending a top element of
$T_2$ to $ux_2 - vx_3$, where

$$\xymatrixcolsep{1.2pc}\xymatrixrowsep{1.6pc}
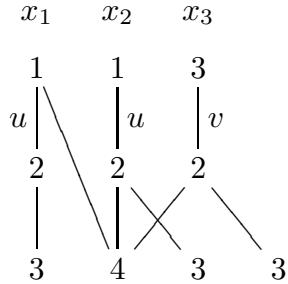
\xymatrix{ 1 \edge[d]_{u} \edge[ddr] {\save+<0ex,4ex> \drop{x_1}
\restore} &1
\edge[d]^{u} {\save+<0ex,4ex> \drop{x_2} \restore} &3 \edge[d]^{v}
{\save+<0ex,4ex> \drop{x_3} \restore} \\ 2 \edge[d] &2 \edge[d]
\edge[dr] &2 \edge[dl] \edge[dr] \\ 3 &4 &3 &3 }$$

\noindent is the graph of $T_4$, relative to suitable top elements
$x_1, x_2, x_3$, normed by $e_1$, $e_1$, and $e_3$, respectively, such
that the socle of
$\la(ux_2 - vx_3)$ equals $S_3$.  One checks that $\nu$ does not belong
to $\etilde_2 \Jtilde^2 \etilde_4$.  On the other hand, the map in
$\Hom_\la(T_2, T_4)$ that sends a top element of
$T_2$ to $ux_1$ factors through $T_5$ by way an obvious choice of  arrow
$\mu: T_2 \rightarrow T_5$ and the map $\delta: T_5 \rightarrow T_4$
specified below.  Consequently, only $\nu$ qualifies as an arrow from
$\etilde_2$ to $\etilde_4$.  We have two linearly independent maps
$\delta$ and $\delta'$ in $\Hom_\la(T_5, T_4)$, both of which belong to
$(\etilde_5 \Jtilde
\etilde_4) \setminus (\etilde_5  \Jtilde^2 \etilde_4)$:  Say
$\delta$ sends a top element $y = e_2 y$ of
$T_5$ to $ux_2$, and $\delta'$ sends $y$ to $vx_3$.  These assignments
pin down $\delta$ and
$\delta'$ up to nonzero scalar factors; in particular, $\delta$ maps
$T_5$ onto a submodule of codimension $1$ of the module $\la x_1 +
\la x_2 \subset T_4$ with graph

$$\xymatrixcolsep{1.0pc}\xymatrixrowsep{1.6pc}
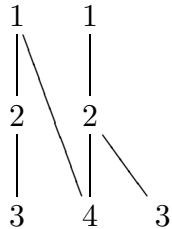
\xymatrix{ 1 \edge[d] \edge[ddr] &1 \edge[d] \\ 2 \edge[d] &2 \edge[d]
\edge[dr] \\ 3 &4 &3 }$$

\noindent and $\delta'$ maps $T_5$ onto a submodule of codimension $1$
of $\la x_1 + \la x_3$. Then the difference $\delta - \delta'$ equals
$\epsilon \nu$ $\in$ $\etilde_5 \Jtilde^2
\etilde_4$, where $\epsilon \in \Hom_\la(T_5, T_2)$ is a suitable
epimorphism sending $y$ to a top element of $T_2$.  So
$\Qtilde$ contains only a single arrow from
$\etilde_5$ to $\etilde_4$; we picked the map $\delta$ for the above
graphs (as opposed to the alternate choice, $\delta'$).
\qed
\enddefinition

\enddocument